\documentclass[10pt,twoside,reqno] {amsart}

\usepackage{amssymb}
\usepackage{pb-diagram} 
\usepackage{graphicx} 

\oddsidemargin=\evensidemargin
\catcode`\@=11
\long\def\@makefntext#1{
\protect\noindent \hbox to 3.2pt {\hskip-.9pt
$^{{\eightrm\@thefnmark}}$\hfil}#1\hfill}       

\def\ps@myheadings{\let\@mkboth\@gobbletwo      
\def\@oddhead{\hbox{}
\rightmark\hfil\eightrm\thepage}
\def\@oddfoot{}\def\@evenhead{\eightrm\thepage\hfil
\leftmark\hbox{}}\def\@evenfoot{}
\def\sectionmark##1{}\def\subsectionmark##1{}}

\catcode`@=11                        
\def\ps@plain{\let\@mkboth\@gobbletwo
     \def\@oddhead{}\def\@oddfoot{\eightrm\hfil\thepage
     \hfil}\def\@evenhead{}\let\@evenfoot\@oddfoot}

\newcounter{sectionc}\newcounter{subsectionc}\newcounter{subsubsectionc}
\renewcommand{\section}[1] {\vspace{12pt}\addtocounter{sectionc}{1}
\setcounter{theorem}{0} \setcounter{equation}{0}
\setcounter{subsectionc}{0}\setcounter{subsubsectionc}{0}\noindent
    {\tenbf\thesectionc. #1}\par\vspace{5pt}}
\renewcommand{\subsection}[1] {\vspace{12pt}\addtocounter{subsectionc}{1}
    \setcounter{subsubsectionc}{0}\noindent
    {\bf\thesectionc.\thesubsectionc.
    {\kern1pt \bfit #1}}\par\vspace{5pt}}
\renewcommand{\subsubsection}[1] {\vspace{12pt}
    \addtocounter{subsubsectionc}{1}
    \noindent
    {\tenrm\thesectionc.\thesubsectionc.\thesubsubsectionc. {\kern1pt
    \it #1}}\par\vspace{5pt}}
\newcommand{\nonumsection}[1] {\vspace{12pt}\noindent{\tenbf #1}
    \par\vspace{5pt}}

\newtheorem{theorem}{Theorem}[sectionc]
\newtheorem{lemma}[theorem]{Lemma}
\newtheorem{proposition}[theorem]{Proposition}
\newtheorem{corollary}[theorem]{Corollary}
\newtheorem{definition}[theorem]{Definition}
\newtheorem{remark}[theorem]{Remark}
\numberwithin{equation}{sectionc}

\newcommand{\publisher}[2]{{\begin{center}\footnotesize\smalllineskip
    Received #1\\
    Revised #2
        \end{center}
    }}

\def\abstracts#1#2#3#4{{
    \centering{\begin{minipage}{4.5in}\footnotesize\baselineskip=10pt
    \centerline{ABSTRACT}
    \parindent=15pt #1\par
    \parindent=15pt #2\par
    \parindent=15pt #3\par
    \parindent=15pt #4\par
    \end{minipage}}\par}}
\def\keywords#1{{
    \centering{\begin{minipage}{4.5in}\footnotesize\baselineskip=10pt
    {\footnotesize\it Keywords}\/: #1
    \end{minipage}}\par}}

\def\subjectclass#1{{
    \centering{\begin{minipage}{4.5in}\footnotesize\baselineskip=10pt
    {\footnotesize\it 2000 Mathematics Subject Classification.}\ #1
    \end{minipage}}\par}}

\newcommand{\textlineskip}{\baselineskip=13pt}
\newcommand{\smalllineskip}{\baselineskip=10pt}

\newcommand{\copyrightheading}[1]
    {\vspace*{-2.5cm}\smalllineskip{\flushleft
    {\footnotesize Journal of Knot Theory and Its Ramifications #1}\\
    {\footnotesize \copyright\kern2pt World Scientific
         Publishing Company}\\
         }}

\renewenvironment{thebibliography}[1]
    {\frenchspacing
     \ninerm\baselineskip=11pt
     \begin{list}{[\arabic{enumi}]}
    {\usecounter{enumi}\setlength{\parsep}{0pt}
     \setlength{\leftmargin 19pt}{\rightmargin 0pt}   
     \setlength{\itemsep}{0pt} \settowidth
    {\labelwidth}{[#1]}\sloppy}}{\end{list}}

\newcommand{\fcaption}[1]{
        \refstepcounter{figure}
        \setbox\@tempboxa = \hbox{\footnotesize Fig.~\thefigure. #1}
        \ifdim \wd\@tempboxa > 5in
           {\begin{center}
        \parbox{5in}{\footnotesize\smalllineskip Fig.~\thefigure. #1}
            \end{center}}
        \else
             {\begin{center}
             {\footnotesize Fig.~\thefigure. #1}
              \end{center}}
        \fi}
\def\runninghead#1#2{\pagestyle{myheadings}
\markboth{{\protect\footnotesize\it{\quad #1}}\hfill}
{\hfill{\protect\footnotesize\it{#2\quad}}}}

\font\tenrm=cmr10

\font\tenbf=cmbx10
\font\bfit=cmbxti10 at 10pt
\font\ninerm=cmr9
\font\nineit=cmti9

\font\eightrm=cmr8

\newcommand{\TryPackage}[3]{\IfFileExists{#1.sty}{\usepackage{#1}#2}{#3}}
\TryPackage{mathrsfs}{\renewcommand{\mathcal}{\mathscr}}{%
     \TryPackage{eucal}{}{}}

\newcommand{\lto}{\longrightarrow}
\newcommand{\al}{\alpha}
\newcommand{\be}{\beta}
\newcommand{\de}{\delta}
\newcommand{\ga}{\gamma}

\newcommand{\varep}{\varepsilon}

\newcommand{\la}{\lambda}

\newcommand{\si}{\sigma}
\newcommand{\Ga}{\Gamma}

\newcommand{\Si}{\Sigma}

\newcommand{\ZZ}{{\mathbb Z}}
\newcommand{\RR}{{\mathbb R}}
\newcommand{\CC}{{\mathbb C}}
\newcommand{\NN}{{\mathbb N}}
\newcommand{\PP}{{\mathbb P}}
\newcommand{\QQ}{{\mathbb Q}}

\newcommand{\cM}{{\mathcal M}}

\newcommand{\mer}{{\mathcal M}} 
\newcommand{\lng}{{\mathcal L}}  
\newcommand{\SLC}{{SL_2({\mathbb C})}}
\newcommand{\slc}{{{\mathfrak s}{\mathfrak l}_2({\mathbb C})}}

\newcommand{\tr}{\operatorname{\it tr}}

\newcommand{\Aut}{\operatorname{Aut}}
\newcommand{\Ad}{\operatorname{\it Ad}}

\date{\today}
\begin{document}
\setlength{\textheight}{7.7truein}  

\runninghead{\quad The SL$_2(\CC)$ Casson Invariant for Seifert fibered homology spheres}
{The SL$_2(\CC)$ Casson Invariant for Seifert fibered homology spheres \quad}

\normalsize\textlineskip
\thispagestyle{empty}
\setcounter{page}{1}

\copyrightheading{}         

\vspace*{0.88truein}

\centerline{\bf THE SL$_2$($\CC$) CASSON INVARIANT FOR SEIFERT FIBERED }
\baselineskip=13pt
\centerline{\bf HOMOLOGY SPHERES AND SURGERIES ON TWIST KNOTS}
\vspace*{0.37truein}
\centerline{\footnotesize HANS U. BODEN}
\baselineskip=12pt
\centerline{\footnotesize\it Department of Mathematics \& Statistics, McMaster University}
\baselineskip=10pt
\centerline{\footnotesize\it Hamilton, Ontario, L8S 4K1 Canada}
\baselineskip=10pt
\centerline{\footnotesize\it {\tt boden@mcmaster.ca}}

\vspace*{10pt}
\centerline{\footnotesize CYNTHIA L. CURTIS}
\baselineskip=12pt
\centerline{\footnotesize\it Department of Mathematics \& Statistics, The College of New Jersey}
\baselineskip=10pt
\centerline{\footnotesize\it Ewing, NJ, 08628 USA}
\baselineskip=10pt
\centerline{\footnotesize\it {\tt ccurtis@tcnj.edu}}

\vspace*{0.225truein}
\publisher{}

\vspace*{10pt}

\vspace*{0.21truein}
\abstracts{We derive a simple closed formula for the $\SLC$ Casson invariant
for Seifert fibered homology 3-spheres
using the correspondence between $\SLC$ character varieties and moduli spaces of parabolic Higgs bundles of rank two. These results are then used to
deduce the invariant for Dehn surgeries on twist knots by
combining computations of the Culler-Shalen
norms with the surgery formula for the $\SLC$ Casson invariant.}{}{}{}

\vspace*{0.21truein}


\smallskip
\subjectclass{Primary: 57M27, Secondary: 57M25, 57M05.}

{}{}{}

\baselineskip=13pt

\section{Introduction}
The 3-manifold invariant  $\la_\SLC(\Si)$  is defined in \cite{C}
by counting certain $\SLC$ representations of $\pi_1 \Si$ and
can be regarded as the $\SLC$ analogue of Casson's invariant.
The goal of this paper is to present two techniques for computing
$\la_\SLC(\Si)$. The first is a direct approach and yields a
simple closed formula for the $\SLC$ Casson invariant for Seifert
fibered homology spheres $\Si(a_1,\ldots,a_n)$. The second
involves the surgery formula of \cite{C} and requires calculation
of the Culler-Shalen seminorms. This approach applies to give
values of the invariant for 3-manifolds obtained by Dehn surgery along
a twist knot.

We begin by introducing some notation and
recalling the definition of the invariant
$\la_\SLC(\Si)$.

Given a finitely generated group $\pi$, denote by $R(\pi)$ the
space of representations $\rho\colon \pi \to \SLC$ and by
$R^*(\pi)$ the subspace of irreducible representations. Recall
from \cite{CS} that $R(\pi)$ has the structure of a complex affine
algebraic set. The {\sl character}  of a representation $\rho$ is
the function $\chi_\rho\colon \pi \to \CC$ defined by setting
$\chi_\rho(\ga)=\tr(\rho(\ga))$ for $\ga \in \pi_1 \Si$. The set
of characters of $\SLC$ representations is denoted $X(\pi)$ and
also admits the structure of a complex affine algebraic set.
Furthermore, there is a canonical projection $t\colon R(\pi) \to
X(\pi)$ defined by $t\colon \rho \mapsto \chi_\rho$ which is
surjective. Let  $X^*(\pi)$ be the subspace of characters of
irreducible representations. Given a manifold $\Si$, we denote by
$R(\Si)$ the variety of $\SLC$ representations of $\pi_1 \Si$ and
by $X(\Si)$ the associated character variety.

Suppose now $\Si$ is a closed, orientable 3-manifold
with a Heegaard splitting $(W_1, W_2, F)$. Here, $F$
is a closed orientable surface embedded in $\Si$, and $W_1$ and $W_2$ are handlebodies
with boundaries $\partial W_1 = F = \partial W_2$
 such that $\Si = W_1\cup_F W_2$.
The inclusion maps $F \hookrightarrow W_i$ and $W_i \hookrightarrow \Si$
induce a diagram of surjections
$$\begin{diagram} 
\node{} \node{\pi_1 W_1}\arrow{se,l}{}\\
\node{\pi_1 F}\arrow{ne,l}{} \arrow{se,l}{} \node{} \node{\pi_1 \Si.}\\
\node{} \node{\pi_1 W_2} \arrow{ne,l}{}
\end{diagram}
$$
This diagram induces a diagram of the associated representation
varieties, where all arrows are reversed and are injective rather
than surjective. On the level of character varieties, this gives a
diagram
 $$\begin{diagram} 
\node{} \node{X(W_1)}\arrow{sw,l}{}\\
\node{X(F)} \node{} \node{X(\Si),}\arrow{nw,l}{} \arrow{sw,l}{} \\
\node{} \node{X(W_2)} \arrow{nw,l}{}
\end{diagram}
$$
of injections.  This identifies $X(\Si)$  as the intersection
$$X(\Si) = X(W_1) \cap X(W_2) \subset X(F).$$

There are natural orientations on all the character varieties
determined by their complex structures. The  invariant
$\la_\SLC(\Si)$ is defined as an oriented intersection number of
$X^*(W_1)$ and $X^*(W_2)$ in $X^*(F)$ which counts only compact,
zero-dimensional components of the intersection. Specifically,
there exist a compact neighborhood $U$ of the zero-dimensional
components of $X^*(W_1)\cap X^*(W_2)$ which is disjoint from the
higher dimensional components of the intersection and an isotopy
$h\colon X^*(F) \to X^*(F)$ supported in $U$ such that
$h(X^*(W_1))$ and $X^*(W_2)$ intersect transversely in $U$. Then
given a zero-dimensional component $\chi$ of the intersection
$h(X^*(W_1))\cap X^*(W_2)$, we may set $\varep_\chi = \pm 1$,
depending on whether the orientation of $h(X^*(W_1))$  followed by
that of $X^*(W_2)$ agrees with or disagrees with the orientation
of $X^*(F)$ at $\chi$.

\begin{definition} Let
$\la_\SLC(\Si) = \sum_\chi \varep_\chi,$
where the sum is  over all zero-dimensional
components $\chi$ of the intersection $h(X^*(W_1))\cap X^*(W_2)$.
\end{definition}

The next result recalls from \cite{C} the basic properties of the
$\SLC$ Casson invariant. Of these, properties (i)--(iv) are
proved explicitly in \cite{C} and (v) is implicit in the definition
of $\la_\SLC(\Si)$. Property (vi),
additivity under connected sum for $\ZZ_2$ homology spheres,
is a consequence of Theorem \ref{connect}.

\begin{theorem}
The invariant $\la_\SLC$ satisfies the following properties:
\begin{enumerate}
\item[(i)] For any 3-manifold $\Si,$ $\la_\SLC(\Si) \geq 0.$
\item[(ii)] If $\Si$ is hyperbolic, then $\la_\SLC(\Si) >0.$
\item[(iii)] If $\la_\SLC(\Si) > 0,$ then there exists an
irreducible representation $\rho\colon \pi_1 \Si \to \SLC$. In
particular, if $\pi_1 \Si$ is abelian, then $\la_\SLC(\Si) = 0.$
\item[(iv)] $\la_\SLC$ satisfies a surgery formula. (See Theorem
\ref{surg-form} for details.)
\item[(v)]  $\la_\SLC(-\Si) =
\la_\SLC(\Si)$, where $-\Si$ is $\Si$ with the opposite
orientation. \item[(vi)] $\la_\SLC(\Si_1 \# \Si_2) =
\la_\SLC(\Si_1) + \la_\SLC(\Si_2)$ for $\ZZ_2$ homology 3-spheres
$\Si_1,\Si_2.$
\end{enumerate}
\end{theorem}

There is an alternative formulation of $\la_\SLC(\Si)$ in terms of
the intersection theory of algebraic varieties \cite{F}. For an
isolated point $\chi$ in $X^*(W_1)\cap X^*(W_2)$, set $m_\chi$
equal to the intersection multiplicity of $\chi$ in the
intersection cycle $X^*(W_1)\cdot X^*(W_2)$.  (See \cite{C} for
details.) Then
\begin{equation} \label{basic}
\la_\SLC(\Si) = \sum_\chi m_\chi,
\end{equation}
where the sum is  over all isolated points
$\chi$ in $X^*(W_1)\cap X^*(W_2)$.

We will use equation (\ref{basic}) to compute $\la_\SLC(\Si)$ in
several interesting cases where the intersection multiplicities
$m_\chi$  are not too difficult to determine.  In Section 2, we
give a simple criterion under which the components of
$X^*(W_1)\cap X^*(W_2)$ are all zero-dimensional with intersection
multiplicity one. In this fortuitous case, one can compute
$\la_\SLC(\Si)$ directly from the character variety $X(\Si)$
without reference to the Heegaard splitting or the isotopy $h$.
After verifying that this criterion holds for Brieskorn spheres,
we use it to calculate $\la_\SLC(\Si),$ first for
Brieskorn spheres, then more generally for Seifert fibered homology
spheres. In Section 3 we state and prove a formula
for  $\la_\SLC$ for connected sums
of rational homology 3-spheres.
In Section 4 we turn our attention to  3-manifolds
resulting from Dehn surgery on a knot in a homology 3-sphere $\Si$
and recall from \cite{C} the surgery formula for the $\SLC$ Casson
invariant. Finally in Section 5 we combine the calculations of
Section 2 with the surgery formula and deduce a formula for the
$\SLC$ Casson invariant for 3-manifolds obtained by Dehn surgery
on a twist knot.

\section{Seifert fibered homology spheres}

Given a representation $\rho\colon \pi_1 \Si\to \SLC$,
we denote by $H^*(\Si;\slc_{\Ad \rho})$  the cohomology groups of  $\Si$
with coefficients in $\slc$ twisted by
$$\Ad  \rho\colon \pi_1 \Si \to \Aut(\slc).$$

\begin{theorem} \label{cohom}
If $\Si$ is a closed oriented 3-manifold and
$\rho\colon \pi_1 \Si \to \SLC$
is an irreducible representation with
$H^1(\Si;\slc_{\Ad \rho})= 0$, then
its character $\chi_{\rho}$ is an
isolated point of the intersection $X^*(W_1)\cap X^*(W_2)$
with intersection multiplicity 1.
\end{theorem}

\noindent
{\bf Proof.}
Let $\rho_0\colon \pi_1 F \lto \pi_1 \Si \stackrel{\rho}{\lto}
\SLC$ be the representation on the Heegaard surface $F$ obtained
by pullback. Similarly, for $i=1,2$ let $\rho_i\colon \pi_1 W_i
\lto \pi_1 \Si \stackrel{\rho}{\lto} \SLC$ be the representation
on the handlebody $W_i$. Because $\rho$ is irreducible and because
each of the three inclusion maps induces a surjection on the level
of fundamental groups, the three representations $\rho_0, \rho_1,$
and $\rho_2$ are irreducible as well.

The Zariski tangent space to $\chi_\rho \in X(\Si)$ is
a subspace of $H^1(\Si; \slc_{\Ad \rho})$. Mayer-Vietoris identifies
this cohomology group with the
intersection of  the images of $H^1(W_1, \slc_{\Ad \rho_1})$ and
$H^1(W_2, \slc_{\Ad \rho_2})$ in
$H^1(F, \slc_{\Ad \rho_0})$.
The condition that $H^1(\Si; \slc_{\Ad \rho})=0$
therefore guarantees that $X(W_1)$ and
$X(W_2)$ intersect transversely at $\chi_\rho$, which is
therefore isolated. Since all spaces are oriented as complex
varieties and the intersection is transverse, this shows that
$\chi_\rho$ contributes $+1$ to $\la_{\SLC}(\Si)$.
\qed\kern0.8pt
\smallskip

\begin{corollary} \label{simple}
If $\Si$ is a closed oriented 3-manifold such that
$H^1(\Si;\slc_{\Ad \rho}) = 0$ for every irreducible
representation $\rho\colon \pi_1 \Si \to \SLC$, then $\la_\SLC(\Si)$ is
an exact count of the conjugacy classes of irreducible
representations of $\pi_1 \Si$ in $\SLC$ -- i.e. $\la_\SLC(\Si) =
|X^*(\Si)|.$
\end{corollary}

We now apply these results to compute
$\la_\SLC(\Si)$ for  Seifert fibered
homology spheres.
We begin with the Brieskorn manifolds
$$\Si(p,q,r)= \{ (x,y,z) \in \CC^3 \mid x^p+y^q+z^r = 0\} \cap S^5.$$
If $p,q,$ and $r$ are positive and pairwise relatively prime, then
$\Si(p,q,r)$ is an integral homology sphere, called a Brieskorn
sphere.

\begin{theorem} \label{pqr}
If $\Si(p,q,r)$ is a Brieskorn sphere and
$\rho\colon \pi_1 \Si(p,q,r) \to \SLC$ is irreducible, then
$H^1(\Si(p,q,r), \slc_{\Ad \rho}) = 0$.
Furthermore, $$\la_\SLC (\Si(p,q,r)) = \frac{(p-1)(q-1)(r-1)}{4}.$$
\end{theorem}
\noindent
{\bf Proof.}
One way to prove this is to use the
correspondence between the character variety $X(\Si)$ and
certain moduli spaces of parabolic Higgs bundles over $\CC\PP^1$.
This is the approach we will adopt in our proof
of Theorem \ref{SFHS}, but here we give an independent and elementary argument.

The fundamental group of $\Si(p,q,r)$ has a presentation
\begin{equation} \label{pqr-pres}
\pi_1 \Si(p,q,r) = \langle x, y, z, h \mid
h \text{ is central, } x^p =h^a, y^q =h^b, z^r =h^c,
xyz=1 \rangle,
\end{equation}
for integers $a,b,c$ satisfying
\begin{equation} \label{eq:abc}
a qr + b pr + c pq = 1.
\end{equation}

\begin{lemma} \label{lemco}
$H^1(\Si(p,q,r);\slc_{\Ad \rho}) =0$ for all $\rho\colon \pi_1 \Si(p,q,r) \to \SLC$.
\end{lemma}
\noindent
{\bf Proof.}
Since $\Si(p,q,r)$ is a homology sphere, the cohomology vanishes
when $\rho$ is the trivial representation.  Thus, we assume $\rho$
is irreducible, and therefore $\rho(h) = \pm I.$ Hence
$$\rho(x)^{2p}= \rho(y)^{2q}=\rho(z)^{2r} = \rho(h)^2=I,$$
so the eigenvalues of $\rho(x), \rho(y),$ and $ \rho(z)$ are
$2p$-th, $2q$-th, and $2r$-th roots of 1, respectively.  Moreover,
since $\rho$ is irreducible, none of these eigenvalues is $\pm 1$.

To simplify notation, we write $\pi$ for $\pi_1 \Si(p,q,r)$. We
use the Fox calculus to determine the dimension of the space of
1-cocycles $ Z^1 (\pi; \slc_{\Ad \rho}).$ The 1-cocycles are
determined by elements $X,Y,Z, H \in \slc$ satisfying the
equations obtained by taking Fox derivatives of the relations in
(\ref{pqr-pres}). The commutation relation $hxh^{-1}x^{-1}$ gives
   $$ H+  \Ad \rho(h)X -  \Ad \rho(x)H-X=0.$$
 Since $ \Ad \rho(h) X = X$,
we see that $H$ lies in the kernel of $1- \Ad \rho(x).$ Similarly,
 $H$ lies in the kernels of  $1- \Ad
\rho(y)$ and $1- \Ad \rho (z).$ Since $\rho$ is irreducible with
image generated
   by $\rho(x),\rho(y)$ and $\rho(z)$, this implies that
$H=0.$

Setting $H=0$ in  the remaining equations, we obtain:
  \begin{eqnarray}
  (1+\Ad \rho(x) + \cdots + \Ad \rho(x^{p-1})) X &=& 0, \label{eq:Adx}\\
  (1+\Ad \rho(y) + \cdots + \Ad \rho(y^{q-1})) Y &=& 0, \label{eq:Ady}\\
  (1+\Ad \rho(z) + \cdots + \Ad \rho(z^{r-1})) Z &=& 0, \label{eq:Adz}\\
  X + \Ad \rho(x) Y + \Ad \rho(xy) Z &=& 0. \label{eq:Adxy}
  \end{eqnarray}

 By (\ref{eq:Adx}), $X$ lies in the kernel of $
(1+\Ad \rho(x) + \cdots + \Ad \rho(x^{p-1}))$.  Now since the
eigenvalues of $\rho(x)$ are distinct $2p$-th roots of 1, we see
that there is an isomorphism $\slc \cong \CC^3$ such that $\Ad
\rho(x)$ is given in the new coordinates by
$$ \Ad \rho(x) (\zeta_1,\zeta_2,\zeta_3) =
(e^{\pi ik/p}\zeta_1,e^{-\pi ik/p}\zeta_2, \zeta_3)$$
for some $1 \leq k < p.$ Then
clearly the kernel of $(1+\Ad \rho(x) + \cdots + \Ad
\rho(x^{p-1}))$ is the 2-dimensional subspace of $\slc$
corresponding to the subspace $\zeta_3=0$ in $\CC^3$. Similarly,
equations (\ref{eq:Ady}) and (\ref{eq:Adz}) imply that the kernels
of $ (1+\Ad \rho(y) + \cdots + \Ad \rho(y^{p-1}))$ and $ (1+\Ad
\rho(z) + \cdots + \Ad \rho(z^{p-1}))$ are 2-dimensional. Finally,
since $\rho$ is irreducible, we see that (\ref{eq:Adxy}) imposes
three independent conditions. We conclude that $$\dim
Z^1(\pi;\slc_{\Ad \rho}) = 2+2+2 -3=3.$$

Now the irreducibility of $\rho$ also implies that
$H^0(\pi;\slc_{\Ad \rho}) =0$. Hence the space of coboundaries has
$\dim B^1(\pi;\slc_{\Ad \rho})=3$, and  it follows that
$$\dim H^1(\pi;\slc_{\Ad \rho})
= \dim Z^1(\pi;\slc_{\Ad \rho}) - \dim B^1(\pi;\slc_{\Ad \rho})
= 0. \hfill \qed\kern0.8pt$$
\smallskip

Now Corollary \ref{simple} applies to show that $\la_\SLC(\Si) =
|X^*(\Si)|$, and we turn to the problem of enumerating the
characters in $X^*(\Si).$ Our enumeration will be given in terms
of the numbers of distinct conjugacy classes of appropriate
non-central roots of unity in $\SLC.$

Reordering $p, q,$ and $r$ as necessary, we may assume that $q$
and $r$ are odd.  We establish a one-to-one correspondence between
characters of irreducible representations $\rho \colon \pi_1
\Si(p,q,r) \to \SLC$ and characters of irreducible representations
$\bar{\rho} \colon T(2p,q,r) \to \SLC$, where $T(2p,q,r)$ is the
triangle group $$T(2p,q,r) = \langle x, y, z \mid x^{2p} = y^q =
z^r = xyz=1 \rangle.$$

Assume first that $\rho \colon \pi_1 \Si(p,q,r) \to \SLC$ is given
and irreducible. If $b$ is odd, we may replace it with $b+q$ and
$a$ with $a-p$ so that $b$ is even. Likewise, we may replace $c$
with $c+r$ and $a$ with $a-p$ as necessary to guarantee that $c$
is even.

Since $b$ is even, $\rho(y)$ must be a $q$-th root of
$\rho(h)^b=I.$ Similarly, since $c$ is even, $\rho(z)$ must be an
$r$-th root of $I$. Now  equation (\ref{eq:abc}) implies $a$ is
odd. Hence $\rho(x)$ is a $p$-th root of $\rho(h)^a =\rho(h)=\pm
I,$ and therefore $\rho(x)$ is a $2p$-th root of $I$. Moreover
$\rho(xyz)=I.$ Thus $\rho$ determines an irreducible
representation $\bar{\rho} \colon T(2p,q,r) \to \SLC$ given by
$\bar{\rho}(x) = \rho(x), \bar{\rho}(y) = \rho(y),$ and
$\bar{\rho}(z) = \rho(z).$

Now assume  $\bar{\rho}\colon T(2p,q,r) \to \SLC$ is given and
irreducible. Setting $\rho(h) = \bar{\rho}(x^p) = \pm I$ defines
an irreducible representation $\rho\colon \pi_1 \Si(p,q,r) \to
\SLC$ in the obvious way.

Thus, $\la_\SLC (\Si(p,q,r)) = |X^*(\Si(p,q,r))| =
|X^*(T(2p,q,r))|.$  We now enumerate the characters in
$|X^*(T(2p,q,r))|.$  Recall that two irreducible representations
$\rho_1, \rho_2:G \to \SLC$ have the same character if and only if
they are conjugate.  Thus, we enumerate the conjugacy classes of
irreducible representations $\rho\colon T(2p,q,r) \to \SLC.$

Note that $T(2p,q,r)$ is generated by the two elements $x$ and
$y,$ since $xyz = 1.$ But for any group $G$ generated by two
elements $g_1$ and $g_2$, an irreducible representation
$\phi\colon G \to \SLC$ is determined up to conjugacy by the
traces of $\phi(g_1), \phi(g_2),$ and $\phi(g_1g_2)$.  To see
this, recall that the Cayley-Hamilton theorem implies the trace
relation:
$$\tr(AB) - \tr(A)\tr(B) + \tr(A^{-1}B) = 0$$
for all $A,B \in \SLC.$ Arguing by induction, this relation implies that
the trace of any word in $A$ and $B$ is determined by
$\tr(A), \tr(B),$ and $\tr(AB)$.

To enumerate the conjugacy classes of irreducible
representations $\rho\colon T(2p,q,r) \to \SLC,$ we need only
determine the possible traces of $\rho(x), \rho(y),$ and $\rho(xy)
= \rho(z^{-1}).$ Since these matrices have order $2p, q$, and $r,$
respectively, in $\SLC$, it is clear that the eigenvalues of these
matrices must be $2p$-th, $q$-th, and $r$-th roots of unity.
Moreover, since $\rho$ is irreducible, none of the eigenvalues are
$\pm 1$.  The next lemma is used to identify the eigenvalues
associated to irreducible representations
$\rho\colon T(2p,q,r) \to \SLC$ .

\begin{lemma} \label{ABC}
Suppose $A,B, C \in \SLC$ satisfy
$\tr A = 2 \cos \al, \tr B = 2 \cos \be,$ and
$\tr C = 2 \cos \ga$ and for $\al,\be,\ga \in (0,\pi)$.
Then
there exists $P \in \SLC$
with $\tr (A P B P^{-1}) = \tr C.$
\end{lemma}
\begin{remark}
It is interesting to compare this result to the $SU(2)$ analogue,
where necessary and sufficient conditions for the existence
of an  $SU(2)$ representation with specified traces
are given by the
quantum Clebsch-Gordon coefficients.
\end{remark}

\noindent
{\bf Proof.}
The trace conditions guarantee that $A,B,$ and $C$ are conjugate
to unitary matrices and as such are diagonalizable; thus we may
assume that
$$A = \begin{pmatrix} e^{i \al} & 0 \\ 0 & e^{-i \al} \end{pmatrix},
\quad B = \begin{pmatrix} e^{i \be} & 0 \\ 0 & e^{-i \be}
\end{pmatrix} \quad \text{and} \quad C = \begin{pmatrix} e^{i \ga}
& 0 \\ 0 & e^{-i \ga} \end{pmatrix}.$$ We will actually find $P
\in SL_2(\RR)$ that satisfies the conclusion by writing $$ P =
\begin{pmatrix} u & v \\ -1 & 1 \end{pmatrix},$$ and computing
\begin{eqnarray*}
\tr (A P B P^{-1}) &=& e^{i(\al+\be)} u + e^{i(\al-\be)} v
+e^{i(\be-\al)} v + e^{-i(\al+\be)} u \\
&=& 2 u \cos(\al+\be) + 2 v \cos(\al-\be).
\end{eqnarray*}
To prove the lemma, we need to find  $u$ and $v$ satisfying
\begin{eqnarray*}
1&=& u + v  \\
\cos \ga &=& u \cos(\al+\be) +  v \cos(\al-\be).
\end{eqnarray*}
Basic linear algebra shows that these equations can be solved unless
$\cos(\al+\be) = \cos(\al-\be),$ which is equivalent to the
condition that $\sin \al = 0$ or $\sin \be=0.$ Since $\al,\be \in
(0,\pi),$ we see that $\sin \al \neq 0 \neq \sin \be.$
\qed\kern0.8pt
\smallskip

Thus, given any 3 elements $\alpha, \beta, \gamma \in (0,\pi)$
such that $e^{i\alpha}$ is a $2p$-th root of unity, $e^{i\beta}$
is a $q$-th root of unity, and $e^{i\gamma}$ is an $r$-th root of
unity, we may define an irreducible representation $\rho\colon T(2p,q,r) \to
\SLC$ with $tr \rho(x) = 2 \cos 2\alpha, tr \rho(y) = 2 \cos
\beta,$ and $tr \rho(xy) = 2 \cos  \gamma$ by setting
$$\quad \rho(x) = \begin{pmatrix} e^{2i \alpha} & 0 \\ 0 & e^{-2i \alpha}
\end{pmatrix} \quad \text{and} \quad \rho(y) = P\begin{pmatrix} e^{i \beta}
& 0 \\ 0 & e^{-i \beta} \end{pmatrix}P^{-1},$$ where $P$ is the
matrix found in Lemma \ref{ABC}.

Finally, noting that there are $p-1$ $2p$-th roots of unity,
$\tfrac{q-1}{2}$ $q$-th roots of unity, and $\tfrac{r-1}{2}$
$r$-th roots of unity in the open semicircle $\{e^{it}\mid
0<t<\pi\},$ we find that
$$\la_\SLC (\Si(p,q,r)) = |X^*(T(2p,q,r))| =
\tfrac{1}{4}(p-1)(q-1)(r-1). \qed\kern0.8pt$$
\smallskip

\begin{theorem} \label{SFHS}
Suppose $a_1,\ldots, a_n$ are positive integers that are
pairwise relatively prime, and denote by  $\Si(a_1, \ldots, a_n)$
the associated Seifert fibered homology sphere. Then
$$\la_\SLC (\Si(a_1,\ldots, a_n)) =\sum_{1 \leq i_1 < i_2 < i_3 \leq n}
 \frac{(a_{i_1}-1)(a_{i_2}-1)(a_{i_3}-1)}{4}.$$
\end{theorem}
\noindent
\begin{remark} The right-hand side of this formula equals $\frac{1}{4}\si_3\,(a_1-1, \ldots, a_n-1)$,
the elementary symmetric polynomial of degree 3
in the $n$ variables $a_1-1, \ldots, a_n-1$.
\end{remark}
\noindent
{\bf Proof.}
We begin with the presentation
\begin{equation} \label{a-pres}
\pi_1 \Si(a_1, \ldots, a_n) =
\langle x_1, \ldots , x_n, h \mid
h \hbox{ central, } x_i^{a_i}=h^{-b_i},
x_1 \cdots x_n = 1\rangle.
\end{equation}
Here, the $b_i$ are not unique but must satisfy
\begin{equation} \label{eq:a}
\sum_{i=1}^n  b_i\, a_1 \cdots \widehat{a}_i \cdots a_n = 1,
\end{equation}
where $\widehat{a}_i$ indicates this term is omitted.

By reordering, we may assume that $a_i$ is odd for $i > 1$. There
is a correspondence between the character variety
$X^*(\Si(a_1,\ldots, a_n))$ and that of the $n$-gon group
$$T(2a_1,a_2,\ldots, a_n) = \langle x_1, \ldots , x_n \mid
x_1^{2a_1}= x_2^{a_2}=\cdots = x_n^{a_n} = x_1 \cdots x_n = 1
\rangle.$$  This correspondence may be established as in the proof
of Theorem \ref{pqr} for the case $n = 3$.

Now Simpson's theorem \cite{S} establishes a correspondence
between connected components of the character variety $X^*(T(2a_1,
a_2, \ldots, a_n))$ and moduli spaces of parabolic Higgs bundles.
Components of $X^*(T(2a_1, a_2, \ldots, a_n))$ are indexed by the
possible traces of the images of the generators $x_1, \ldots,
x_n$. Since each of them has finite order, any irreducible
representation must map each $x_i$ to a unitary matrix with
eigenvalues given by roots of unity of the appropriate order.
Thus, the traces of the images of $x_1, x_2, \ldots, x_n$ are real
numbers of the form $$2\cos(2 \pi \al_1), 2\cos(2 \pi \al_2),
\ldots, 2\cos(2 \pi \al_n),$$ where $\al_1 = k_1/2 a_1$ and $\al_i =
k_i/a_i$ for $i >1,$ and where $k_i \in \ZZ$ satisfies $0 \leq k_1
\leq a_1$ for $i=1$ and $0 \leq k_i < a_i/2$ for $i>1.$

We set $\al = (\al_1,\ldots, \al_n)$ and let $\cM_\al$ denote the
moduli space of rank two parabolic Higgs bundles of parabolic
degree zero over $\CC\PP^1$ with $n$ marked points $p_1, \ldots,
p_n$ and weights $\al_i,1-\al_i$ at $p_i$. Given this choice of
weights, one can easily verify that every semistable parabolic
Higgs bundle is actually stable. As a result, the moduli space
$\cM_\al$ is smooth of complex dimension $2m - 6$, where
$$m = m(\al)=\left| \left\{ \al_i \mid \al_i \in \left(0,\tfrac{1}{2}\right)\right\} \right|,$$
the number of nontrivial flags in the quasiparabolic structure.
(If $m < 3$, then  $\cM_\al = \varnothing$.
This corresponds to the requirement that an irreducible
representation of $T(2a_1,a_2,\ldots, a_n)$ must send
at least three generators to noncentral elements in $\SLC$.)

Thus, the zero-dimensional components of $X^*(T(2a_1,a_2,\ldots,
a_n))$ are in one-to-one correspondence with the subset $\{ \al
\mid m(\al)=3\}$ of all possible weights. Such a weight
$\al=(\al_1,\ldots, \al_n)$ is obtained by choosing $1\leq i_1 <
i_2 < i_3 \leq n$ with $\al_{i_1}, \al_{i_2}, \al_{i_3} \in
\left(0,\tfrac{1}{2}\right)$ and all remaining $\al_j$ equal to
$0$ or $\tfrac{1}{2}$. (Note: only $\al_1$ is allowed to equal
$\tfrac{1}{2}$ here, since for $\chi \in X^*(T(2a_1,a_2,\ldots,
a_n))$, the only generator whose trace may equal $-2$ is $x_1$.
Note also that as in the proof of Theorem \ref{pqr}, we may apply Lemma \ref{ABC}
to show that there exists an irreducible representation of $T(2a_1,a_2, \ldots, a_n)$
with the given weights.)

A straightforward generalization of Lemma \ref{lemco} shows that
such components correspond to isolated components of
$X^*(\Si(a_1,\ldots, a_n))$ of intersection multiplicity one, and
so our task is simply to enumerate them. To do this, consider the
cases $i_1 =1$ and $i_1 >1$ separately.

In the case $i_1 = 1,$ we have
$$\al_1 \in \left\{ \tfrac{1}{2a_1}, \ldots, \tfrac{a_1-1}{2a_1} \right\},$$
a set with $a_1-1$ elements. (This enumeration corresponds to the
one used previously in counting conjugacy classes of roots of
unity in $\SLC$. In particular, there are $a_1-1$ distinct
conjugacy classes of noncentral $2a_1$-th roots of unity in
$\SLC.$) Likewise,
$$\al_{i_2} \in \left\{ \tfrac{1}{a_{i_2}}, \ldots, \tfrac{a_{i_2}-1}{2a_{i_2}} \right\}
\quad \text{and} \quad
\al_{i_3} \in \left\{ \tfrac{1}{a_{i_3}}, \ldots, \tfrac{a_{i_3}-1}{2a_{i_3}} \right\},$$
which are sets with $\frac{1}{2}(a_2-1)$ and $\frac{1}{2}(a_3-1)$ elements,
respectively. (Note that $a_{i_2}$ and $a_{i_3}$ are both odd,
and that  $\frac{1}{2}(a_2-1)$ and $\frac{1}{2}(a_3-1)$
are precisely the number of  distinct conjugacy classes
of noncentral $a_{i_2}$-th
and $a_{i_3}$-th roots of unity in $\SLC$, respectively.)
Since $m(\al)=3,$ we take all other $a_j=0,$ so
this case give a  total of
$\frac{1}{4}(a_{i_1}-1)(a_{i_2}-1)(a_{i_3}-1)$
isolated point components.

In the case $i_1>1,$ one counts as above to see that for $k=1,2,3,$
$$\al_{i_k} \in \left\{ \tfrac{1}{a_{i_k}}, \ldots,
\tfrac{a_{i_k}-1}{2a_{i_k}}\right\},$$ a set with
$\frac{1}{2}(a_{i_k}-1)$ elements. In this case, $\al_1$ may equal
0 or $\tfrac{1}{2},$ which gives an extra factor of 2, and all
other $\al_j=0.$  Thus, there are
$\frac{1}{4}(a_{i_1}-1)(a_{i_2}-1)(a_{i_3}-1)$ isolated point
components. The proof of the theorem is completed by summing over
all possible $1 \leq i_1 < i_2 < i_3 \leq n.$
\qed\kern0.8pt
\smallskip

For Brieskorn spheres, Fintushel and Stern showed that Casson's
invariant satisfies $$\la_{SU(2)}(\Si(p,q,r) =  \tfrac{1}{8}{\rm
signature}\, M(p,q,r),$$ where $M(p,q,r)$ denotes the Milnor fiber of
the singularity \cite{FS}. This equation was shown to hold more
generally by Neumann and Wahl, who have conjectured that
$\la_{SU(2)}(\Si) = \frac{1}{8}{\rm signature}\, M$ for any link of a
normal complete intersection singularity \cite{NW}.

Theorem \ref{pqr}
shows that  the $\SLC$ Casson invariant satisfies
\begin{equation} \label{milnor}
\la_\SLC(\Si(p,q,r)) = \tfrac{1}{4} \mu(p,q,r),
\end{equation}
where $\mu(p,q,r)$ denotes the Milnor number \cite{M}. Theorem
\ref{SFHS} shows this formula does not extend to Seifert fibered
homology spheres $\Si(a_1,\ldots, a_n)$ with $n \geq 4$. We explain
briefly the correct generalization of (\ref{milnor}).

Consider the components $X_j \subset X^*(\Si(a_1,\ldots, a_n))$ of
highest dimension, namely those with $\dim X_j = 2n-6$.  By the
formula on p.597 of \cite{BY}, we have $\chi(X_j) = (n-1)(n-2)
2^{n-4}$ for each component $X_j$.  Moreover, a lattice point
count similar to the argument given in the proofs of Theorems
\ref{pqr} and \ref{SFHS} shows that there are $\tfrac{1}{4}(a_1 -
1)\ldots(a_n - 1)$ such components.  Finally, by Theorem 9.1 of
\cite{M} we have $\mu(a_1,\ldots,a_n) = (a_1 - 1)\ldots(a_n - 1).$
It follows that
$$\sum_j \chi(X_j) = (n-1)(n-2) 2^{n-6} \mu(a_1, \ldots, a_n).$$

\section{Connected sum formula}
\begin{theorem} \label{connect}
If $\Si_1$ and $\Si_2$ are rational homology spheres,
then
$$\la_\SLC(\Si_1 \# \Si_2) = |H_1(\Si_2; \ZZ_2)| \cdot \la_\SLC(\Si_1) +
|H_1(\Si_1; \ZZ_2)| \cdot  \la_\SLC(\Si_2),$$
 As a consequence,
the invariant $\la_\SLC$ is additive under connected sum for
$\ZZ_2$-homology spheres.
\end{theorem}
\noindent
{\bf Proof.}
Since $\pi_1 (\Si_1 \# \Si_2) = \pi_1 (\Si_1)* \pi_1(\Si_2)$, we
have $R(\Si_1 \# \Si_2) = R(\Si_1) \times R(\Si_2).$ Suppose
$\rho_1, \rho_2$ are $\SLC$ representations on $\Si_1, \Si_2$,
with stabilizer groups $\Ga_1, \Ga_2$, respectively.  Then under
the conjugation action, the orbit of $\rho_1$ in $R(\Si_1)$ is
$\SLC/\Ga_1$ and the orbit of $\rho_2$  in $R(\Si_2)$ is
$\SLC/\Ga_2$. We are interested only in $\SLC$ representations
$\rho = (\rho_1, \rho_2) \in R(\Si_1 \# \Si_2)$ which are
irreducible and map to isolated points under $t\colon R(\Si_1 \#
\Si_2) \to X(\Si_1 \# \Si_2).$

Given a representation $\rho = (\rho_1,\rho_2) \in R(\Si_1 \#
\Si_2)$ and an element $g\in \SLC,$ we may define a representation
$\rho_g \in R(\Si_1 \# \Si_2)$ by  $\rho_g = (g \cdot \rho_1,
\rho_2)$. In general, $\rho_g$ lies in a different conjugacy class
from $\rho$. In fact, the conjugacy class $[\rho]$ lies on a
family $[\rho_g]$ of conjugacy classes parameterized by the double
coset $\Ga_1 \backslash \SLC / \Ga_2,$ and the conjugacy classes
of irreducible representations in this family are in one-to-one
correspondence with characters. In particular, if an irreducible
representation $\rho = (\rho_1,\rho_2)$ maps to an isolated point
in $X(\Si_1 \# \Si_2)$, then either $\Ga_1$ or $\Ga_2$ must equal
$\SLC$. This is the case precisely when either $\rho_1$ or
$\rho_2$ is central.

Suppose then that $\rho_2$ is central, so the image of $\rho_2$ is
in $\{ \pm I \}$. It follows easily that $\rho = (\rho_1,\rho_2)$
is irreducible if and only if $\rho_1$ is irreducible. Further,
the resulting character $\chi_\rho$ will be isolated in $X^*(\Si_1
\# \Si_2)$ if and only if $\chi_{\rho_1}$ is isolated in
$X^*(\Si_1).$ (Here, we use the Mayer-Vietoris theorem and the
fact that $H^1(\Si_2;\slc_{\Ad \rho_2})=0$, which follows since
$\rho_2$ is central and $\Si_2$ is a rational homology 3-sphere.)
One can piece together a compactly supported perturbation from
$\Si_1$ with the trivial perturbation on $\Si_2$ to make the
intersection transverse in a neighborhood of $\chi_\rho,$ showing
in essence that $m_{\chi_\rho} = m_{\chi_{\rho_1}}$.

An analogous argument shows that if $\rho_1$ is central, then
$\rho = (\rho_1,\rho_2)$ is irreducible and isolated if and only
if $\rho_2$ is irreducible and isolated, and if so then
$m_{\chi_\rho} = m_{\chi_{\rho_2}}$.  The theorem follows from the
additional observation that the central representations of $\pi_1
\Si_i$ are in one-to-one correspondence with elements of
$H_1(\Si_i; \ZZ_2).$
\qed\kern0.8pt

\newpage

\section{Dehn surgery formula for small knots}
Theorem 4.8 of \cite{C} gives a formula for the $\SLC$ Casson
invariant for Dehn surgeries on small knots in homology spheres in
terms of a weighted sum of Culler-Shalen seminorms. In this
section, we restate that formula, incorporating the comments from
\cite{C1}. We begin by reviewing the notation of \cite{CS, CCGLS}.

Suppose $M$ is a compact, irreducible, orientable 3-manifold with
boundary a torus. An {\sl incompressible surface} in $M$ is a
properly embedded surface $(F,\partial F) \to (M,\partial M)$ such
that $\pi_1 F \to \pi_1 M$ is injective and no component of $F$ is
a 2-sphere bounding a 3-ball. An {\sl essential surface} in $M$ is
an incompressible surface in $M,$ no component of which is
boundary parallel. The manifold $M$ is called {\sl small} if it
does not contain a closed essential surface, and a knot $K$ in
$\Si$ is called {\sl small} if its complement $\Si \setminus
\tau(K)$ is a small manifold.

If $\ga$ is a simple closed curve in $\partial M$, the Dehn
filling of $M$ along $\ga$ will be denoted by $M(\ga)$; it is the
closed 3-manifold obtained by identifying a solid torus with $M$
along their boundaries so that $\ga$ bounds a disk. Note that the
homeomorphism type of $M(\ga)$ depends only on the {\sl slope} of
$\ga$ -- that is, the unoriented isotopy class of $\ga$. Primitive
elements in $H_1(\partial M; \ZZ)$ determine slopes
under a two-to-one correspondence.

If $F$ is an essential surface in $M$
with nonempty boundary, then all of its  boundary
components are parallel and the slope of one (hence all) of these
curves is called the {\sl boundary slope} of $F$. A
slope is called a {\sl strict boundary slope} if it is
the boundary slope of an essential surface that is not the
fiber of any fibration of $M$ over $S^1$.

For each $\ga \in \pi_1 M,$ there is a regular map $I_\ga\colon X(M)
\to \CC$ defined by $I_\ga(\chi) = \chi(\ga).$ Let
$e\colon H_1(\partial M;\ZZ)\to \pi_1(\partial M)$ be the
inverse of the Hurewicz
isomorphism.  Identifying $e(\al) \in \pi_1(\partial M)$ with its
image in $\pi_1 M$ under the natural map $\pi_1(\partial M) \to
\pi_1 M,$ we obtain a well-defined function $I_{e(\al)}$ on $X(M)$
for each $\al \in H_1(\partial M; \ZZ).$ Let $f_\al\colon X(M) \to \CC$
be the regular function defined by $f_\al  = I_{e(\al)} -2$ for
$\al \in H_1(\partial M; \ZZ).$

\begin{remark}
Our function $f_\al$ does not agree with that of the papers
\cite{CS, CGLS, BZ1, BZ2}, where $f_\al$ denotes $I^2_{e(\al)} -
4.$
\end{remark}

Let $r\colon X(M)\to X(\partial M)$ be the restriction map induced by
$\pi_1(\partial M) \to \pi_1 M$. Suppose $X_i$ is a
one-dimensional component of $X(M)$ such that $r(X_i)$ is
one-dimensional.  Let $f_{i,\al}\colon X_i \to \CC$ denote the regular
function obtained by restricting $f_{\al}$ to $X_i$ for each $i$.

Let $\widetilde{X}_i$ denote the smooth, projective curve
birationally equivalent to $X_i$. Regular functions on $X_i$
extend to rational functions on $\widetilde{X}_i$. We abuse
notation and denote the extension of $f_{i,\al}$ to
$\widetilde{X}_i$ by $f_{i, \al}\colon \widetilde{X}_i \to \CC
\cup \{ \infty \} = \CC \PP^1.$

A generalization of the argument in Section 1.4 of \cite{CGLS}
proves the following result.

\begin{proposition} Suppose $X_i$ is a one-dimensional component of
$X(M)$ which contains an irreducible character and whose
restriction $r(X_i)$ is also one-dimensional. There exists a
seminorm  $\| \cdot \|_i$ on the real vector space $H_1(\partial
M; \RR)$ satisfying
$$ \| \al \|_i = \tfrac{1}{2}
\deg(f_{i,\al})$$
for all $\al$ in the lattice  $H_1(\partial M; \ZZ)$.
\end{proposition}

\begin{remark}
The relationship between $\|   \cdot   \|_i$ and the Culler-Shalen
norm introduced in \cite{CGLS} is as follows:  Suppose $M$ is
hyperbolic. Let $X_0$ be the unique component  of $X(M)$
containing the character of the $\SLC$ lift of the associated
discrete, faithful representation $\rho\colon \pi_1 M \to PSL_2(\CC)$.
Then Culler, Gordon, Luecke and Shalen proved the existence of a
norm $\| \cdot \|_{CS}$ on $H_1(\partial M; \RR)$ such that $\|
\al \|_{CS} = \deg \left( I_{0,e(\al)}^2 - 4\right)$ for all $\al
\in H_1(\partial M; \ZZ).$ (See Section 1.4 of \cite{CGLS}.) It
follows immediately from the definitions that $\| \cdot \|_{CS}= 4
\| \cdot \|_0.$ A similar relationship holds between $\| \cdot
\|_i$ and the Culler-Shalen seminorms described in \cite{BZ1,
BZ2} associated to the remaining curves $X_i$.
\end{remark}

We will relate the $\SLC$-Casson invariant of a manifold obtained
by Dehn filling to this seminorm; however we must impose certain
restrictions on the filling slope.

\begin{definition} The slope of a simple closed curve $\ga$ in $\partial M$
is called {\sl irregular} if there exists an irreducible representation
$\rho\colon \pi_1 M \to \SLC$ such
that
\begin{itemize}
 \item [(i)] the character $\chi_{\rho}$ of $\rho$ lies on a one-dimensional
component $X_i$ of $X(M)$ such that $r(X_i)$ is one-dimensional,
 \item [(ii)] $\tr \rho (\al) =\pm 2$ for all
 $\al$ in the image of $i^*\colon \pi_1(\partial M) \to \pi_1(M),$
 \item[(iii)] $\ker (\rho \circ i^*)$ is the cyclic group generated by $[\ga] \in \pi_1(\partial M)$.
\end{itemize}
A slope is called {\sl regular} if it is not irregular.
\end{definition}

\begin{remark} In \cite{C}, irregular slopes
are called {\sl exceptional} slopes.
However, the term ``exceptional slope" is
commonly used to refer to a slope $\ga$ along
which the Dehn filled
manifold $M(\ga)$ is not hyperbolic.
The term irregular is used here
to avoid any ambiguity.
\end{remark}

With these definitions, we are almost ready to state the
Dehn surgery formula. But first we recall some
useful notation for Dehn fillings along knot complements. For any
choice of basis $(\al,\be)$  for $H_1(\partial M; \ZZ),$ there is
a bijective correspondence between unoriented isotopy classes of
simple closed curves in $\partial M$ and elements in $\QQ \cup \{
\infty\}$ given by $\ga\mapsto p/q,$ where  $\ga = p \al + q \be$.
If $M$ is the complement of a knot $K$ in an integral homology
sphere $\Si$, then the meridian $\mer$ and longitude $\lng$ of $K$
provide a preferred basis for $H_1(\partial M; \ZZ).$ We set
$K(p/q)=M(\ga),$  the Dehn filling along the curve $\ga = p \mer +
q \lng.$ In this case, we call $p/q$ the slope of $\ga$ and
$K(p/q)$ the result of $p/q$ Dehn surgery along the knot $K$ in
the homology 3-sphere $\Si$.

\begin{definition}
A slope $p/q$ is called {\sl admissible} for $K$ if
\begin{enumerate}
\item[(i)] $p/q$ is a regular slope which is not a strict boundary slope.
\item[(ii)] no $p'$-th root of unity is a root of the Alexander polynomial of
$K$, where $p'=p$ if $p$ is odd and $p' = p/2$ if $p$ is even.
\end{enumerate}
\end{definition}

The next result is a restatement of Theorem 4.8 of \cite{C}, as
corrected in \cite{C1}.
\begin{theorem} \label{surg-form}
Suppose $K$ is a small knot in a homology 3-sphere $\Si$ with
complement $M$. Let $\{X_i\}$ be the collection of all
one-dimensional components of the character variety $X(M)$ such
that $r(X_i)$ is one-dimensional and such that $X_i \cap X^*(M)$
is nonempty. Define $\si\colon \ZZ \to \{0,1\}$ by $\si(p) =0$ if
$p$ is even and $\si(p)=1$ if $p$ is odd.

Then there exist integral weights $n_i > 0$ depending only on
$X_i$ and non-negative numbers $E_0$ and $E_1$ in $\frac{1}{2}
\ZZ$ depending only on $K$ such that for every admissible slope
$p/q$, we have

$$
\la_\SLC (K(p/q) ) = \sum_i n_i \| p\mer + q \lng \|_i - E_{\si(p)}.
$$
\end{theorem}


In fact each character $\chi \in X^*(W_1)\cap X^*(W_2)$
contributes equally to each side of the equation in Theorem
\ref{surg-form}.

\begin{proposition} \label{int-mult}
Suppose $K$ is a small knot in a homology 3-sphere $\Si$ with
complement $M$. Let $\al = p \mer + q \lng \in H_1(M;\ZZ)$ and
suppose $\chi\in X^*(M)$ extends to give an isolated character in
$X^*(K(p/q))$. Assume further that $\chi(\mer) \neq \pm 2$ or
$\chi(\lng) \neq \pm 2$. Then the intersection multiplicity of
$X^*(W_1)$ and $X^*(W_2)$ at $\chi$ is given by
$$m_\chi =\tfrac{1}{2} \sum_i n_i \mu_{i,\chi},$$
where  $\mu_{i,\chi}$ is the order of vanishing of
$f_{i,\al}(\chi)=0$. In particular,
it follows that
 $$m_\chi \geq \sum_{\{i \mid \chi \in X_i\}} n_i.$$
\end{proposition}

A proof of Proposition \ref{int-mult} is included
in Section 6 for completeness.

\medskip

Let $\pi^*\colon X(K(p/q)) \to X(M)$ be the injective map induced by the
surjection $\pi\colon \pi_1 M \to \pi_1(K(p/q))$.

\begin{corollary} \label{cor-int-mult}
Let $K$ be a small knot in an integral homology sphere with
complement $M$.  Suppose $\rho\colon \pi_1(K(p/q)) \to \SLC$ is an
irreducible representation such that $H^1(K(p/q); \slc_{\Ad
\rho})=0$.  Let $\chi_{\rho}$ be the character of $\rho$, and
suppose either $\chi_{\rho}(\mer) \neq \pm 2$ or
$\chi_{\rho}(\lng) \neq \pm 2$. Then there is a unique curve $X_i$
in $X(M)$ containing $\pi^* \chi_{\rho}$, and the
 associated weight of $X_i$ is $n_i = 1$.
\end{corollary}
\noindent
{\bf Proof.}
 Since $H^1(K(p/q); \slc_{\Ad \rho})=0$, we know
that the intersection multiplicity of $X^*(W_1)$ and $X^*(W_2)$ at
$\chi_{\rho}$ is 1 by Theorem 2.1.  Then the assertion follows
immediately from the proposition.
\qed\kern0.8pt
\smallskip

Thus, finding a character $\chi_\rho \in X^*(K(p/q))$ with
$H^1(K(p/q);\slc_{\Ad \rho})=0$ and $\chi_{\rho}(\mer) \neq \pm 2$
or $\chi_{\rho}(\lng) \neq \pm 2$ determines the weight $n_i$ for
the corresponding curve $X_i \subset X(M)$.  This weight in turn
yields information about the intersection multiplicities for all
characters coming from $X_i$ in the closed manifolds $K(a/b)$ for
all regular, non-boundary slopes $a/b$. In particular, if every
curve $X_i$ in $X(M)$ contains a character which is the image of
the character of
 a representation $\rho$ with
$\chi_{\rho} \in X^*(K(p/q))$ for some $p/q$ and
with $H^1(K(p/q);\slc_{\Ad \rho})=0$, then one can determine
$\la_\SLC(K(a/b))$ for all admissible slopes $a/b$
directly from information about the
seminorms $\| \cdot \|_i$ and the numbers $E_0, E_1$.  We
apply this in the next section to manifolds obtained by Dehn
surgery on twist knots.

\section{Dehn surgeries on twist knots}

\begin{figure}[h]
\begin{center}
\leavevmode\hbox{}
\includegraphics[height=1.3in]{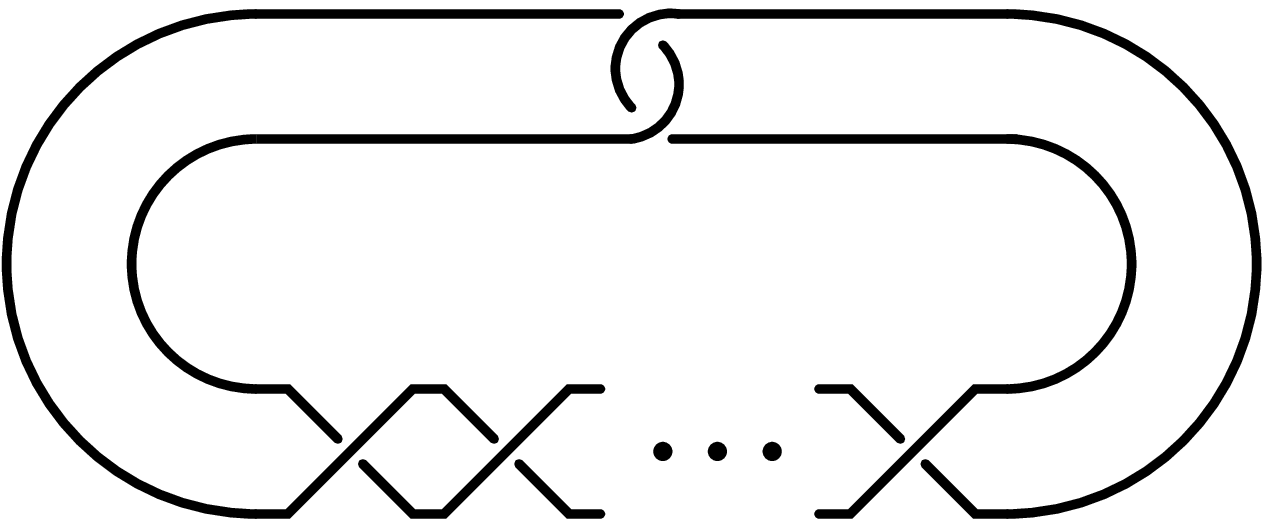}
\fcaption{{The twist knot $K_\xi$ with a clasp and $\xi$ half twists.}} \label{ttwist}
\end{center}
\end{figure}

Let $K_\xi$ for $\xi \geq 1$ be the $\xi$-twist knot in $S^3$
depicted in Figure \ref{ttwist}, with complement $M_\xi$, and let
$K_\xi(p/q)$ denote the 3-manifold obtained by performing $p/q$
Dehn surgery on $K_\xi$. In this section, we present  formulas for
$\la_\SLC(K_\xi(p/q))$ for all slopes $p/q$ which are not strict
boundary slopes. Note that $K_\xi$ is hyperbolic for $\xi
> 1.$

By \cite{Bu}, all irreducible characters in $X(M_\xi)$ lie on a
single curve in $X(M_\xi)$ which we denote $X_\xi$. We denote by
$\| \cdot \|_\xi$ the seminorm on $H_1(\partial M_\xi; \RR)$
obtained from this curve. (Note that when $\xi > 1$, $\| \cdot
\|_\xi$ is actually a norm, since $K_\xi$ is hyperbolic.) Let
$n_\xi$ denote the weight of $X_\xi$ described in Theorem
\ref{surg-form}.

\begin{proposition}\label{twistsurg}
Suppose $p/q$ is an admissible slope for $K_\xi$.  Then
$$\la_\SLC(K_\xi(p/q)) = \| p\mer + q\lng \|_{\xi} - E_{\si(p)}.$$
\end{proposition}

\noindent
{\bf Proof.}   By Theorem \ref{surg-form} above, since $X(M_\xi)$ contains a
single curve $X_\xi$ containing irreducible characters, we must
show that $n_{\xi} = 1$.

It is well-known that $K_\xi$ is small. In fact, if $\xi = 2k-1$
is odd, then
$$K_{2k-1}(1) = \Si(2,3,6k - 1).$$
One can demonstrate this by identifying both manifolds with the
one obtained by $1/k$ Dehn surgery on the left-handed trefoil
using Kirby calculus. Similarly, if $\xi = 2k$ is even, then
$$K_{2k}(-1) = -\Si(2,3,6k + 1).$$

By Theorem \ref{pqr}, each of these Seifert fibered spaces
$\Sigma$ admits irreducible representations $\rho\colon \pi_1 \Si
\to \SLC$, and each such $\rho$ satisfies $H^1(\Sigma;\slc_{\Ad
\rho}) = 0$.  We show that each of these Seifert fibered spaces
admits an irreducible representation $\rho$ with character
$\chi_{\rho} \in X^*(\Sigma)$ satisfying $\chi_{\rho}(\mer) \neq
\pm 2$. Then Corollary \ref{cor-int-mult} will imply that $n_{\xi}
= 1$.

We use Casson's $SU(2)$ invariant to establish the existence of
such characters. We have
$$-\la_{SU(2)}(K_{2k-1}(1)) = k = \la_{SU(2)}(K_{2k}(-1)).$$
Thus there exists an irreducible $SU(2)$ representation $\rho$ of
$\pi_1(\Sigma)$, since $k\neq 0$. But $\pi_1(\Sigma)$ is normally
generated by the meridian $\mer$, so $\rho(\mer) \neq \pm I$ since
$\rho$ is irreducible.  Since $\rho(\mer)\in SU(2)$, it follows
that $\chi_{\rho} (\mer) \neq \pm 2$.
\qed\kern0.8pt
\smallskip

Thus, $\lambda_{\SLC}(K_{\xi}(p/q))$ is completely determined by
the seminorm $\| \cdot \|_{\xi}$ and the numbers $E_0$ and $E_1$
for admissible slopes. These will be computed in Propositions
\ref{CS-twist} and \ref{E0E1}, but first, we investigate
admissibility for slopes of twist knots.

The first lemma establishes regularity for all slopes of twist knots.
\begin{lemma} \label{reg}
If $K_\xi$ is a nontrivial twist knot, then all
slopes are regular.
\end{lemma}
\noindent
{\bf Proof.}
We have the presentation $\pi_1(M_\xi)=\langle x,y \mid
xw_{\xi}=w_{\xi}y\rangle$, where $w_{\xi}
=(yx^{-1}y^{-1}x)^{\xi/2}$ if $\xi$ is even, and $w_{\xi} =
(yxy^{-1}x^{-1})^{(\xi-1)/2} yx$ if $\xi$ is odd. Then as in
section 7 of \cite{CCGLS},  the curve $X_\xi$ in $X(M_\xi)$ can be
parameterized by characters of representations $\rho\colon
\pi_1(M_\xi)\to \SLC$ with
$$\rho(x) = \left[ \begin{array}{cc}
\mu & 1\\
0 & \mu^{-1} \end{array} \right] \quad \text{and} \quad \rho(y) =
\left[ \begin{array}{cc}
\mu & 0\\
t_{\xi} & \mu^{-1} \end{array} \right], $$ where $t_{\xi}$ is
chosen so that $\rho(xw_{\xi}) = \rho(w_{\xi}y)$.  Henceforth we
will write $w$ for $w_{\xi}$ and $t$ for $t_{\xi}$ whenever no
confusion will arise by our doing so.

With this parameterization of $X_\xi$, we may identify the
meridian $\mer$ with $x$, so $\rho(\mer) = \rho(x)$. If $w^*$
denotes the word obtained by reversing $w$, then the longitude
$\lng$ is given by $\lng = ww^*$ if $\xi$ is even and by
$\lng=x^{-4}ww^*$ if $\xi$ is odd.

Suppose $p/q$ is an irregular slope.  Then there exists an
irreducible representation $\rho$ of the knot complement with
$\chi_{\rho}(\mer) = \pm 2, \chi_{\rho}(\lng) = \pm 2,$ and $
\rho(\mer^p \lng^q) = I.$  We show that this cannot occur.

Since $\chi(\mer) =\pm 2$, we know that $\mu = \mu^{-1}= \pm 1$.
Suppose
$$\rho(w)=\left[\begin{array}{cc} \al & \be \\ \ga & \de \end{array} \right].$$
Then
$$\rho(xw)=\left[ \begin{array}{cc} \mu & 1 \\ 0 & \mu \end{array} \right]
\left[ \begin{array}{cc} \al & \be \\ \ga & \de
\end{array} \right] =
\left[ \begin{array}{cc} \mu \al +\ga & \mu \be +\de \\ \mu \ga &
\mu \de
\end{array} \right]$$
and
$$ \rho(wy)=
\left[ \begin{array}{cc} \al & \be \\ \ga & \de
\end{array} \right] \left[ \begin{array}{cc} \mu & 0 \\ t & \mu
\end{array} \right]
=\left[ \begin{array}{cc}
\mu\al + t \be & \mu \be \\
\mu \ga + t \de &  \mu \de
\end{array} \right].$$

Note that $\rho(xw)=\rho(wy)$ implies that $\de=0$. It follows
that $\ga = - \be^{-1}$ since $\rho(w) \in \SLC$.

Now the presentation of $\pi_1(M_{\xi})$ implies that $w_{\xi + 2}
= y x^{-1} y^{-1} x w_{\xi}$ if $\xi$ is even and $w_{\xi + 2} = y
x y^{-1} x^{-1} w_{\xi}$ if $\xi$ is odd.  Then simple proofs by
induction on $\xi$ (treating $\xi$ odd and $\xi$ even separately)
show that
 $$\rho(w^*)= \left[\begin{array}{cc}0 & \be
\\  -\be^{-1} & \al
\end{array}\right].$$  Hence
$$
 \rho(ww^*)   =  \left[\begin{array}{cc}\al & \be \\
 -\be^{-1} & 0 \end{array}\right]\left[\begin{array}{cc} 0 & \be \\
 -\be^{-1} & \al \end{array}\right]
    =  \left[\begin{array}{cc} -1 & 2 \al \be \\ 0 & -1
   \end{array} \right]
   $$
and
$$ \rho(x^{-4} ww^*)   =  \left[\begin{array}{cc}\mu & -1 \\
 0 & \mu \end{array} \right]^4
 \left[\begin{array}{cc} -1 & 2 \al \be \\ 0 & -1
   \end{array} \right]
=\left[\begin{array}{cc} -1 & 2 \al \be + 4 \mu^3 \\ 0 & -1
   \end{array} \right].$$
Thus we have
$$\rho(\mer)=\left[\begin{array}{cc} \pm1 & 1 \\ 0 & \pm 1
   \end{array} \right] \quad \text{and} \quad \rho(\lng)=\left[\begin{array}{cc} -1 & u \\ 0 & -1
   \end{array} \right]$$
 where $u=2\alpha \beta$ if $\xi$ is even and $u = 2 \alpha \beta + 4 \mu^3$ if $\xi$ is odd.
  An easy computation shows that
 if $\chi(\mer) =2$, then
\begin{equation} \label{evpq}
\rho(\mer^p \lng ^q) = \left[\begin{array}{cc} 1 & 1 \\ 0 & 1
   \end{array} \right]^p \left[\begin{array}{cc} -1 & u \\ 0 & -1
   \end{array} \right]^q = (-1)^q \left[\begin{array}{cc} 1 & p-qu \\ 0 &  1
   \end{array} \right].
   \end{equation}
 On the other hand,
   if $\chi(\mer) =-2,$ then
\begin{equation} \label{odpq}
\rho(\mer^p \lng ^q) =  \left[\begin{array}{cc} -1 & 1 \\ 0 & -1
   \end{array} \right]^p \left[\begin{array}{cc} -1 & u \\ 0 & -1
   \end{array} \right]^q = (-1)^{p+q} \left[\begin{array}{cc} 1 & -p-qu \\ 0 &  1
   \end{array} \right].
   \end{equation}

\noindent {\bf Claim.} If $u$ is rational, then it is an even
integer.

\medskip
\noindent Assume $u \in \QQ$.
Using that $u = 2\alpha \beta$
when $\xi$ is even
and that $u = 2 \alpha \beta + 4 \mu^3$
when $\xi$ is odd, together
with the fact that $\mu = \pm 1,$
the assumption implies that $\al \be \in \QQ$.
The claim will follow once we show $\al \be \in \ZZ.$

Arguing by induction on $\xi$ as before, it is not difficult to
establish that for each $\xi \in \NN$ and for $\mu = \pm 1$:
\begin{enumerate}
\item[(i)] $\al(t), \be(t), \ga(t)$, and $\de(t)$ are all
polynomials with integer coefficients. \item[(ii)]$\de(t)$ is a
monic polynomial. \item[(iii)] $\ga(t) = t \be(t)$. \item[(iv)] If
$\xi$ is even, then $\al(t) - \de(t) = (2 - t)\be(t)$. \item[(v)]
If $\xi$ is odd, then $\al(t) + \de(t) = (2+t)\be(t).$
\end{enumerate}

Fix $\xi$. Recall that $t_{\xi}$ is chosen so that
$\rho(x w_{\xi}) = \rho(w_{\xi} y)$, which implies that
$\de(t_xi)=0$.
Assuming that
$t_{\xi} \in \QQ$, then it follows
from the fact that $t_{\xi}$ is the root of a monic
polynomial with integer coefficients that $t_\xi \in \ZZ$, and \
then (i) implies that $\al(t_{\xi})\be(t_{\xi}) \in \ZZ$
as claimed.
Thus the claim will follow once we show that $t_\xi \in \QQ.$

Since $\ga(t_{\xi}) = -\be(t_{\xi})^{-1}$, (iii) implies that
$t_{\xi} \beta(t_{\xi})^2 = -1$. Also
$\al(t_{\xi}) = \al(t_{\xi})-\de(t_{\xi}) =
\al(t_{\xi})+\de(t_{\xi})$. Now if $\xi$ is even, then (iv) shows
that
$$\al(t_{\xi})\be(t_{\xi})=(\al(t_{\xi})-\de(t_{\xi}))\be(t_{\xi})=
(2-t_{\xi})\beta(t_{\xi})^2
=2\be(t_{\xi})^2+1.$$
Similarly, if $\xi$ is odd, then (v) gives that
$$\al(t_{\xi})\be(t_{\xi})=(\al(t_{\xi})+\de(t_{\xi}))\be(t_{\xi})=
(2+t_{\xi})\be(t_{\xi})^2
=2\be(t_{\xi})^2-1.$$
In either case, we can solve for $\be(t_{\xi})^2$ in terms of
$\al(t_{\xi}) \be(t_{\xi})$. Since $\al(t_{\xi})\be(t_{\xi}) \in
\QQ$, we have $\be(t_{\xi})^2 \in \QQ.$
The equation $t_{\xi} \be(t_{\xi})^2 = -1$ then
implies that $t_{\xi} \in \QQ$, and
this completes the proof of the claim.
\medskip

By the claim,
either $u \not \in \QQ$ or $u$ is an even integer.
Assume first that $\chi(\mer) =2.$
If $u \not \in \QQ$, then $p \neq qu$ and $\rho(\mer^p \lng^q) \neq I$
by equation (\ref{evpq}).  On the other hand, if $u$
is an even integer and $\rho(\mer^p \lng^q) = I$,  then
equation (\ref{evpq}) implies that $q$ is even
and that $p=qu$.
It follows that $p$ is also even, which
contradicts the assumption that $p$ and $q$ are relatively
prime.

Assume now that $\chi(\mer) =-2.$
If $u \not \in \QQ$, then $p \neq -qu$ and $\rho(\mer^p \lng^q) \neq I$
by equation (\ref{odpq}).
On the other hand, if $u$
is an even integer and $\rho(\mer^p \lng^q) = I$,
 then equation (\ref{odpq}) implies that $p+q$ is even
and that $p=-qu$. Hence $p$ is even, as is $q$, again
contradicting the assumption that $p$ and $q$ are relatively
prime. This completes the proof of the lemma.
\qed\kern0.8pt
\smallskip

The next lemma discusses the roots of
Alexander polynomials of twist knots.

\begin{lemma} \label{alex}
The Alexander polynomial $\Delta_{K_\xi}(t)$ of $K_\xi$ has a
root which is a $k$-th root of unity if and only if $\xi = 1$.
The roots of $\Delta_{K_1}(t)$ are 6-th roots of unity.
\end{lemma}
\noindent
{\bf Proof.}
The Alexander polynomial of $K_\xi$ is given by
$$2 \Delta_{K_\xi}(t) = \begin{cases}
 \xi t^2 - 2(\xi+1)t + \xi &\text{if $\xi$ is even,} \\
 (\xi+1)t^2 - 2\xi t + (\xi+1) &\text{if $\xi$ is odd.}
\end{cases}$$

If $\xi=1$, then $\Delta_{K_1}(t) = t^2 - t + 1$ with roots
$e^{\pm {\pi i}/3}$, which are 6-th roots of unity. We consider the
remaining cases separately.

\noindent {\bf Case 1:} $\xi>1$ is even.

\noindent
Then the roots of $\Delta_{K_\xi}(t)$ are real and equal to $({\xi+1 \pm
\sqrt{2\xi+1}})/\xi$, neither of which has absolute value equal to 1.

\noindent {\bf Case 2:}  $\xi>1$ is odd.

\noindent
Then the roots of $\Delta_{K_\xi}(t)$ are not real. Suppose
$\Delta_{K_\xi}(t)$ has a root of the form $\al = e^{{2 \pi i j}/k}$. Then
$2 \Delta_{K_\xi}(t) = (\xi+1) (t - \al) (t - \bar{\al})$.
Thus, $2\Delta_{K_\xi}(t)$ divides the polynomial $(\xi+1) (t^k -
1)$. It follows that $2\Delta_{K_\xi}(t)/(\xi+1) = t^2 - [2\xi/(\xi+1)]t +
1$ divides $t^k - 1$. We will show this leads to a contradiction.

 Let $t^2 - \be t + 1$ be a quadratic polynomial which divides
$t^k - 1$. We claim that if $\be$ is rational, then it is
integral. Proving this claim will complete the proof of Case 2,
since ${2\xi}/({\xi+1})$ is rational but not integral for $\xi >
1$.

To verify the claim, we write $t^k - 1 = (t^2 - \be t + 1) \cdot
Q(t)$ and solve for the coefficients of $Q(t) = \sum_{i = 0}^{k -
2} q_i t^{k - i - 2}$. We find that $q_0 = 1$, $q_1 = \be$, and
$q_i = \be q_{i-1} - q_{i - 2}$ for $i \geq 2$.  Thus $q_i =
q_i(\be)$ is a monic polynomial of degree $i$ in $\be$ with
integer coefficients.  Moreover $q_{k-2} = -1$, so $\be$ is a root
of the monic polynomial $q_{k-2}(x) + 1$.  But any rational root
of a monic polynomial with integer coefficients is integral. This
establishes the claim and thereby the lemma.
\qed\kern0.8pt
\smallskip

These lemmas allow us to precisely identify the admissible slopes
for $K_{\xi}$:

\begin{proposition} \label{admiss}
Every slope  is admissible for $K_1$ except the strict boundary
slope and slopes of the form $p/q$ where $p$ is a nonzero multiple
of 12. For $\xi>1,$ every slope is admissible for $K_\xi$ except
the strict boundary slopes.

The only strict boundary slope for $K_1$ is $6$.

The strict boundary slopes for $K_2$ are $-4$ and $4$.

The strict boundary slopes for $K_\xi, \xi \ge 3$ odd, are  $0,4,$
and $2\xi +4.$

The strict boundary slopes for $K_\xi, \xi \ge 4$ even, are
$-4,0,$ and  $2\xi .$
\end{proposition}
\noindent
{\bf Proof.}
The first two claims follow from Lemmas \ref{reg} and \ref{alex},
so all that remains is to verify the lists of strict boundary
slopes. The boundary slopes for $K_{\xi}$ are determined in
\cite{HT} to be 0 and 6 if $\xi = 1$; $-4, 0,$ and $2\xi$ if $\xi
\ge 2$ is even; and $0, 4,$ and $2\xi + 4$ if $\xi\ge 3$ is odd.
But the trefoil $K_1$ and the figure eight knot $K_2$ are fibered
knots with 0 the slope of the fiber.  Thus, the only strict
boundary slope for the trefoil is 6, and the strict boundary
slopes of the knots $K_{\xi}$ for $\xi > 1$ are as given above.
\qed\kern0.8pt
\smallskip

In light of Proposition \ref{twistsurg}, we find that computing
the seminorm $\| \cdot \|_{\xi}$ and $E_\si, \ \si =0,1,$ for
$K_\xi$ completely determines $\la_\SLC(K_\xi(p/q))$ for all
admissible slopes.

We compute the Culler-Shalen norm $\| \cdot \|_{CS}$ for $K_\xi,
\, \xi >1$, which determines  $\|\cdot\|_{\xi}$ since
$\|\cdot\|_{\xi}=(1/4)\|\cdot\|_{CS}$. These norms were computed
previously by Boyer, Mattman, Zhang in \cite{BMZ}.  An alternative
method for computing these norms was developed by Ohtsuki in
\cite{O}.

\begin{proposition} \label{CS-twist}
The Culler-Shalen norm for $K_\xi, \, \xi > 1,$ is given by
$$ \| p\mer + q\lng\|_{CS} =
\begin{cases}
 \xi |4q+p| + (\xi-2)|p| + 2|2\xi q - p|  &\text{if $\xi$ is even,} \\
  (\xi-1)|4q-p| +(\xi-1)|p| + 2|(2\xi +4)q-p| &\text{if $\xi$ is odd.}
  \end{cases}$$
\end{proposition}

\noindent
{\bf Proof.}  Recall that $K_{\xi}$ is hyperbolic, so $\|\cdot\|$ is a norm on $X(M_{\xi})$.
In the proof of Proposition 8.6 of \cite{BZ2}, it is shown that
$$\|p\mer+q\lng\|_{CS} = \sum a_i |u_iq - v_i p|,$$
where the sum is taken over all boundary slopes $u_i/v_i$ and
where $a_i \in 2\ZZ$ is non-negative.

Again using the enumeration of the boundary slopes for $K_\xi$ in
\cite{HT}, we have

$$\| p\mer + q\lng\|_{CS} = \begin{cases}
 a_1 |4q+p| + a_2|p| + a_3|2\xi q - p|  &\text{if $\xi$ is even,}\\
 b_1|4q-p| +b_2|p| + b_3|(2\xi+4)q-p| &\text{if $\xi$  is odd.}
\end{cases}
$$

In \cite{BMZ}, the authors compute $\| \cdot \|_{CS}$ for the
slopes $\mer$, $-\mer+\lng$, and $-2\mer+\lng$ for each $K_\xi$
with $\xi$ even and for the mirror image of each $K_\xi$ with
$\xi$ odd. Recalling that $p/q$ Dehn surgery on a knot is
equivalent to $-p/q$ Dehn surgery on its mirror image, we may use
these values to determine $\| \cdot \|_{CS}$ as follows:

 If $\xi$ is even we obtain the equations
\begin{eqnarray*}
2\xi & = & \|\mer\|_{CS} = a_1 +  a_2 +  a_3  \\
8\xi& = & \| - \mer + \lng\|_{CS} = 3 a_1 +  a_2 + (2\xi+1)a_3   \\
8\xi & = & \| - 2\mer + \lng\|_{CS} = 2 a_1 + 2 a_2 + (2\xi+2) a_3
\end{eqnarray*}
with solution
$a_1 = \xi$, $a_2 = \xi-2$, and $a_3 = 2$
easily determined by linear algebra.

If $\xi$ is odd, we obtain the equations
\begin{eqnarray*}
2\xi & = & \|-\mer\|_{CS}=  b_1 +  b_2 +  b_3\\
8\xi+2 & = &  \| \mer + \lng\|_{CS}= 3 b_1 +  b_2 + (2\xi+3)b_3 \\
8\xi & = &  \| 2\mer + \lng\|_{CS}=  2 b_1 + 2 b_2 + (2\xi+2)b_3
\end{eqnarray*}
with solution $b_1 = \xi-1=b_2$ and $b_3 = 2$.
\qed\kern0.8pt
\smallskip

We next compute the correction terms $E_0$ and $E_1$ for the twist
knots.
\begin{proposition} \label{E0E1}
For the twist knots $K_\xi, \ \xi \geq 1,$
we have $E_0 = 0$ and  $E_1 = \xi/2$.
\end{proposition}
\noindent
{\bf Proof.}
We begin with $E_1$.
Note that the
meridian $\mer$ is a regular slope and
is not a boundary slope for $K_\xi$ for any $\xi$.
Therefore by Theorem \ref{surg-form} we have
\begin{eqnarray*}
\la_\SLC(K_\xi(1/0)) & = & \|\mer\|_{\xi} - E_1 \\
& = & 1/4 \|\mer\|_{CS} - E_1.
\end{eqnarray*}
But $K_\xi(1/0)$ is simply $S^3$, so $E_1 = 1/4 \|{\mer}\|_{CS}=\xi/2$
by Proposition \ref{CS-twist}.

We now turn to the computation of $E_0$. By Theorem
\ref{surg-form}, we know that
$$E_0 = \| 2\mer + \lng\|_{\xi} - \la_\SLC(K_\xi(2/1)).$$
Moreover, by Proposition \ref{int-mult} and Corollary
\ref{cor-int-mult}, we know that all characters $\chi$ with
$\chi(\mer) \neq \pm 2$ or $\chi(\lng) \neq \pm 2$ contribute
equally to both $\|2\mer + \lng\|_{\xi}$ and $\la_\SLC(K_\xi(2/1))$.
Thus, $E_0$ measures the difference in the contributions of
characters $\chi$ with $\chi(\mer) =\pm 2$ and $\chi(\lng) =\pm 2$
to the quantities $\|2\mer+\lng\|_{\xi}$ and
$\la_\SLC(K_\xi(2/1))$.  However, by the proof of Lemma \ref{reg},
for any character $\chi$ with $\chi(\mer) =\pm 2$ we have
$\chi(2\mer+\lng) = -2$. Thus, there are no characters $\chi$ with
$\chi(\mer) =\pm 2$ and $\chi(\lng) =\pm 2$ satisfying
$f_{2\mer+\lng}(\chi)=0.$ By Proposition 4.2 of \cite{C}, it
follows that there are no characters $\chi$ with $\chi(\mer) =\pm
2$ and $\chi(\lng) = \pm 2$ which contribute to
$\la_\SLC(K_\xi(2/1))$. Therefore all characters $\chi$ with
$\chi(\mer) = \pm 2$ and $\chi(\lng) = \pm 2$ contribute 0 to both
$\|2\mer+\lng\|_{\xi}$ and $\la_\SLC(K_\xi(2/1))$, and $E_0 = 0$.
\qed\kern0.8pt
\smallskip

Summarizing, we have the following formulas for
$\la_\SLC(K_\xi(p/q))$:

\begin{theorem} \label{twist}
Suppose $\xi > 1$ and $p/q$ is not a strict boundary slope for $K_{\xi}$.  Then the
$\SLC$ Casson invariant of $K_\xi(p/q)$ is as follows:

If $\xi$ is even, then
$$\la_\SLC(K_\xi(p/q))=
\begin{cases}
\frac{1}{4} (\xi |4q+p| + (\xi-2)|p| + 2|2\xi q - p| ) &
  \text{if $p$ is even,}\\
\frac{1}{4} (\xi |4q+p| + (\xi-2)|p| + 2|2\xi q - p| - 2\xi)  &
   \text{if $p$ is odd.}
   \end{cases}$$

If $\xi$ is odd, then
$$\la_\SLC(K_\xi(p/q))=
\begin{cases}
\frac{1}{4} \left((\xi-1)(|4q-p| +|p|) + 2|2\xi q+4q-p| \right)
    &\text{if $p$ is even,}\\
\frac{1}{4} \left((\xi-1)(|4q-p| +|p|) + 2|2\xi q+4q-p| -2 \xi\right) &
\text{if $p$ is odd.}
\end{cases}$$
\end{theorem}

The next two results treat the trefoil $K_1$ as a special case.

\begin{proposition} \label{CS-trefoil}
 The Culler-Shalen norm
for the right hand trefoil $K_1$ is given by
$$ \| p\mer + q\lng\|_{CS} =
 2 |6q-p|$$
\end{proposition}
\noindent
{\bf Proof.}
Proposition \ref{admiss} shows that 6 is the only strict boundary
slope for $K_1$. Using the parameterization of $X^*(M_1)$
described in the proof of Lemma \ref{reg}, it is easy to check
that $\chi(\mer^6 \lng) = -2$ for every $\chi \in X^*(M_1)$.
Therefore by Proposition 5.4 of \cite{BZ1}, $\| \cdot \|_{CS}$ is
an indefinite seminorm which can be written as $\|p \mer + q \lng
\|_{CS}= m|6q-p|$ for some $m \in \ZZ$. Theorem \ref{pqr} gives
that $\la_\SLC(K_1(1)) = \la_\SLC(\Si(2,3,5)) = 2$, and
Proposition \ref{E0E1} implies $E_1 = 1/2.$ Hence $\|\mer +
\lng\|_{1} = 5/2$ and $\|\mer + \lng\|_{CS} =10.$ Thus $m=2$, and
$\| \cdot \|_{CS}$ is as stated in the proposition.
\qed\kern0.8pt
\smallskip

\begin{theorem} \label{trefoil}
The $\SLC$-Casson invariant for the right hand trefoil $K_1$ is
given by
$$\la_\SLC(K_1(p/q)) =
\begin{cases}
\frac{1}{2}|6q-p| - \frac{1}{2} & \text{if $p$ is odd,}\\
\frac{1}{2}|6q-p| & \text{if $p$ is even and not a multiple of $12$,} \\
\frac{1}{2}|6q-p| - 2 & \text{if $p$ is a multiple of $12$.}
\end{cases}$$
\end{theorem}
\noindent
{\bf Proof.}
The theorem follows from Propositions \ref{twistsurg}, \ref{E0E1},
and \ref{CS-trefoil}, provided $p$ is not a multiple of 12. Of
course, one must check the strict boundary slope  $p/q = 6$ by
hand. Since $K_1(6) = L(2,1) \# L(3,1)$, its fundamental group
$\pi_1(K_1(6))=\ZZ_2 * \ZZ_3$ does not admit irreducible $\SLC$
representations and we see that
$\la_\SLC(K_1(6))=0$, in agreement with our formula.
In what follows, we present the argument in
the case $p \neq 0$ is a multiple of 12.

By Proposition \ref{int-mult}, for
any character $\chi$ in $X^*(K(p/q))$ such that either $\chi(\mer)
\neq \pm 2$ or $\chi(\lng) \neq \pm 2$, the intersection multiplicity of
$X^*(W_1)$ and $X^*(W_2)$ at $\chi$ is 1/2 the order of vanishing
of $f_{p\mer+q\lng}$ at $\chi$, since $n_1 = 1$.

Now parameterizing  $X_1$ as in the proof of Lemma \ref{reg}, we
see by the argument used in that proof that for any representation
$\rho\colon \pi_1(M)\to \SLC$ with character $\chi_{\rho}$ in
$X^*(K_1(p/q))$ such that $\chi_{\rho}(\mer)=\pm 2$, the matrix
$\rho(\mer^p)$ is upper triangular with diagonal entries equal to
1, and the matrix $\rho(\lng^q)$ is upper triangular with diagonal
entries equal to $-1$, since $p$ is a multiple of 12 and $q$ is
odd.  Thus, no characters $\chi$ with $\chi(\mer) =\pm 2$ and
$\chi(\lng) =\pm 2$ contribute to $\la_\SLC(K_1(p/q)).$ Therefore
$\la_\SLC(K_1(p/q))$ is 1/2 the sum of the orders of the zeros of
$f_{p \mer +q\lng}$ at characters $\chi$ in $X^*(K_1(p/q))$.

We next observe that $X_1$ contains the characters of exactly two
conjugacy classes of reducible representations: namely, the
representations $\rho_1$ and $\rho_2$ where
 $$\rho_1(x) = \left[ \begin{array}{cc}
e^{\pi i/6} & 1\\
0 & e^{-\pi i/6} \end{array} \right] \quad \text{and} \quad
\rho_1(y) = \left[ \begin{array}{cc}
e^{\pi i/6} & 0\\
0 & e^{-\pi i/6} \end{array} \right]$$
and
$$ \rho_2(x) = \left[ \begin{array}{cc}
e^{5\pi i/6} & 1\\
0 & e^{-5\pi i/6} \end{array} \right] \quad \text{and} \quad
 \rho_2(y) = \left[ \begin{array}{cc}
e^{5\pi i/6} & 0\\
0 & e^{-5\pi i/6} \end{array} \right]. $$  The order of vanishing
of $f_{p\mer+q\lng}$ at each of the associated characters is 2
(which can be seen using an argument analogous to that used to
prove Proposition \ref{int-mult}).  Therefore
$$\la_\SLC(K_1(p/q))= \tfrac{1}{2}\deg( f_{p \mer+q\lng}) - 2 =
\tfrac{1}{4} \| p\mer + q\lng\|_{CS} - 2 = \tfrac{1}{2}|6q - p| -
2.$$
\qed\kern0.8pt
\smallskip

Theorems \ref{pqr} and \ref{twist}
can be used to deduce the knot invariant $\la'_{\SLC}(K)$
defined in \cite{C} for torus and twist knots.

\begin{corollary}
(i) If $T_{p,q}$ is the $p,q$ torus knot,
then $$\la'_{\SLC}(T_{p,q}) = \tfrac{1}{4}pq(p-1)(q-1).$$
(ii) If $K_\xi$ is the $\xi$-twist knot,
then $$\la'_{\SLC}(K_\xi) = \begin{cases}
2 \xi & \text{if $\xi$ is even,}\\
2 \xi+1 & \text{if $\xi$ is odd.}
\end{cases}$$
\end{corollary}
\noindent
{\bf Proof.}
For any small knot $K$ in an integral homology sphere, the
invariant $\la'_{\SLC}(K)$ is defined by the equation
$$ \la'_{\SLC}(K) =
\la_{\SLC}\left(K\left(\tfrac{1}{n+1}\right)\right) -
\la_{\SLC}\left(K\left(\tfrac{1}{n}\right)\right)$$ for $n$ large.
To prove part (i), we use the well-known fact that $1/n$ Dehn
surgery on $T_{p,q}$ gives the Brieskorn sphere $\Si(p,q,pqn-1),$
and Theorem \ref{pqr} implies that
$$ \la_{\SLC}(\Si(p,q,pqn-1)) = \tfrac{1}{4}(p-1)(q-1)(pqn-2).$$
Part (ii) follows similarly from Theorems \ref{twist} and \ref{trefoil}.
\qed\kern0.8pt
\smallskip

It is an important problem to compute $\la_{\SLC}(K(p/q))$
for surgeries on a broader collection of knots, say for
2-bridge knots and/or alternating knots. This
data will help determine the relationship between
$\la'_{\SLC}(K)$ and classical knot
invariants such as the
Alexander polynomial and signature.

\section{Proof of Proposition \ref{int-mult}} \label{app}
This section is
devoted to the proof of Proposition \ref{int-mult}. Essentially, we review
the proof of Theorem 4.8 of \cite{C}, explaining the correction of \cite{C1}.

We begin with some definitions.

Let $\Delta \subset R (\partial M)$ be the subvariety of diagonal
representations, and let $t_\Delta\colon \Delta \to  X (\partial M)$
be the restriction to $\Delta$ of the canonical surjection
$t\colon R(\partial M) \to X (\partial M)$. Then  $t_{\Delta}$ is easily
seen to be surjective. Let $r\colon X(M) \to X (\partial M)$ be the
restriction map. Finally, for each component $X_i$ of $X (M) $
such that $r(X_i)$ is one--dimensional, let $D_i$ be the curve
$t_\Delta^{-1} (\overline{r(X_i)}) \subset \Delta$.

Note that $D_i$ is a branched double cover of the closure of
$r(X_i)$. We have
$$
 \overline{r(X_i)} = t_{\Delta} (D_i)= \frac{D_i}{(a,b)
\sim (a^{-1} , b^{-1})} $$ where $(a,b) \in \CC^*  \times
\CC^* \cong \Delta $.

\bigskip
\noindent {\em Proof of Proposition \ref{int-mult}.} Note that by
Propositions 4.1 and 4.2 of \cite{C}, since $\chi(\mer) \neq \pm
2$ or $\chi(\lng) \neq \pm 2$, we have $\chi \in X^*(W_1)\cap
X^*(W_2)$ if and only if $\chi$ is a zero of
$f_{p\mer+q\lng}.$ By Proposition 4.3 of
\cite{C}, the contribution of $r(\chi)$ to $\la_\SLC(K(p/q))$ is
$$\sum_{\{i \mid \chi \in X_i\}} n_i d_i \phi_{i,\chi},$$ where $n_i$ is a
positive integer depending only on $X_i$, where $d_i$ is the
degree of the map $r|_{X_i}\colon X_i \to r(X_i)$, and where
$\phi_{i,\chi}$ is the intersection number of $r(X_i)$ and
$t_{\Delta}(\{(x,y) \mid x^p y^q = 1\})$ at $r(\chi)$.
Equivalently, $\phi_{i,\chi}$ is 1/2 the intersection multiplicity
of $D_i$ and the curve $\{(x,y) \mid x^p y^q = 1\}$ in $\CC^*
\times \CC ^*$, since $r(\chi)$ has two inverse images in
$t^{-1}_{\Delta}(r(X_i))=D_i$. Thus, the contribution of $\chi$ to
$\la_\SLC(K(p/q))$ is
\[\sum_{\{X_i \mid \chi \in X_i\}} \frac{ n_i d_i \phi_{i,\chi}}
{|r^{-1}(r(\chi))|},\]
where $|r^{-1}(r(\chi))|$ denotes the number of elements in the set
$r^{-1}(r(\chi))$.

 Our
task is to show that
$\mu_{i,\chi} = \frac{2d_i \phi_{i,\chi}}{|r^{-1}(r(\chi))|}$.

For convenience, set $\al = p \mer + q \lng.$
Let $g_{i,\al}\colon r(X_i)\to \CC$ be the map taking an
element $r(\chi) \in r(X_i)$ to $\chi(\mer^p \lng^q)$. Then
$f_{i,\al}$ factors as $g_{i,\al} \circ r$.
Thus,
$$\mu_{i,\chi} = \frac{d_i}{|r^{-1}(r(\chi))|}
\cdot \left(\text{ order of vanishing of $g_{i,\al}$
at $r(\chi)$}\right).$$
Therefore, to prove the first assertion of the proposition, we
must show that the order of the vanishing of $g_{i,\al}$
at $r(\chi)$ is the intersection multiplicity
of $D_i$ and the curve $\{(x,y) \mid x^p y^q = 1\}$ in
$\CC^* \times \CC^*$.

This argument is a pointwise version of Proposition 6.6, Corollary
6.7, and Proposition 8.5 of \cite{BZ2}.

Consider the commutative diagram:
$$\begin{diagram} 
\node{D_i \subset \CC^* \times \CC^*} \arrow{s,l}{\Psi_{p,q}}
\arrow{e,t}{t_\Delta}
\node{r(X_i)} \arrow{s,r}{g_{i,\al}}\\
\node{\CC^*} \arrow{e,t}{z \mapsto z+z^{-1}-2} \node{\CC}
\end{diagram}
$$

Here,  $\Psi_{p,q}\colon \CC^* \times \CC^* \to \CC^*$ is the map
$(x,y) \mapsto x^p y^q.$ Observe that the two horizontal maps in
this diagram have degree two.

We see that $r(\chi)$ is a zero of $g_{i,p\mer+q\lng}$ if
and only if $(x,y) \in t^{-1}_{\Delta}(r(\chi))$ satisfies
$\Psi_{p,q}(x,y)= 1$, and if so, then the orders of vanishing coincide.
But the solutions to $\Psi_{p,q}(x,y)=1$ are precisely the points in the
intersection $D_i \cap \{(x,y) \mid x^p y^q = 1\}$, and the order
of vanishing equals the intersection multiplicity
of $D_i$ and $\{(x,y) \mid x^p y^q = 1\}$ in $\CC^* \times \CC^*$.

This proves the first assertion in the proposition. The final
statement follows from the fact that $\mu_{i,\chi}/2$ is integral.
To see this, note that $\frac{d_i}{|r^{-1}(r(\chi))|}$ is
integral, and since
 $r(\chi)$ has two preimages in $D_i$, it follows that
$\phi_{i,\chi}$ integral. \qed\kern0.8pt

\nonumsection{Acknowledgements}
We would like to thank Steve Boyer, Andr\'{a}s
N\'{e}methi, and Andy Nicas for several illuminating discussions.
The first author was supported by the Natural Sciences and Engineering Research Council of Canada and is grateful to the Max Planck Institute for Mathematics for its hospitality.


\nonumsection{References}


\begin{thebibliography}{999}

\bibitem{BY}
H. U. Boden and K. Yokogawa, {\nineit Moduli spaces of parabolic Higgs
bundles and parabolic $K(D)$ pairs over smooth curves: I}, Intern.
J. Math.  {\bf 7} (1996), 573--598.

\bibitem{BMZ}
S. Boyer, T. Mattman, and X. Zhang, {\nineit The fundamental polygons
of twist knots and the $(-2,3,7)$ pretzel knot,} Knots '96, World
Scientific Publishing Co. (1997), 159--172.

\bibitem{BN}
S. Boyer and A. Nicas,
{\nineit Varieties of group representations and
Casson's invariant for rational homology 3-spheres,}
Trans. Amer. Math. Soc. {\bf 322} (1990), 507--522.

\bibitem{BZ1}
S. Boyer and X. Zhang,
{\nineit On Culler-Shalen seminorms and Dehn filling,}
Annals of Math. {\bf 148} (1998), 737--801.

\bibitem{BZ2}
S. Boyer and X. Zhang,
{\nineit A proof of the finite filling conjecture,}
J. Diff. Geom. {\bf 59} (2001), 87--176.

\bibitem{BW}
M. Brittenham and Y.-Q. Wu,
{\nineit The classification of exceptional Dehn surgeries on 2-bridge knots,}
Comm. Anal. Geom. {\bf 9} (2001), 97--113.

\bibitem{Bu}
G. Burde,
{\nineit  $SU(2)$ representation spaces for two-bridge knot groups,}
Math. Ann. {\bf 288} (1990), 103--119.

\bibitem{CCGLS}
D. Cooper, M. Culler, H. Gillet, D. D. Long, and P. B. Shalen,
{\nineit  Plane curves associated to character varieties of 3-manifolds,}
Invent. Math. {\bf 118} (1994),   47--84.

\bibitem{CS}
M. Culler and P. B. Shalen,
{\nineit Varieties of group representations and splittings of 3-manifolds,}
Annals of Math. {\bf 117} (1983), 109--146.

\bibitem{CGLS}
M. Culler, C. McA. Gordon, J. Luecke, and P. B. Shalen,
{\nineit  Dehn surgery on knots,}
Annals of Math. {\bf 125} (1987), 237--300.

\bibitem{C}
C. L. Curtis, {\nineit  An intersection theory count of the
$\SLC$-representations of the fundamental group of a 3-manifold,}
Topology {\bf 40} (2001), 773--787.

\bibitem{C1}
C. L. Curtis, {\nineit Erratum to ``An intersection theory count of the
$\SLC$-representations of the fundamental group of a 3-manifold,"}
Topology {\bf 42} (2003), 929.

\bibitem{FS}
R. Fintushel and R. Stern,
{\nineit  Instanton homology of Seifert fibred homology three spheres,}
Proc. London Math. Soc. {\bf 61} (1990), 109--137.

\bibitem{F}
W. Fulton, Intersection Theory, Springer-Verlag, Berlin Heidelberg
(1984).

\bibitem{HT}
A.  Hatcher and W. Thurston,
{\nineit  Incompressible surfaces in 2-bridge knot complements,}
Invent. Math. {\bf 79} (1985), 225--246.

\bibitem{M}
J. Milnor,
Singular Points of Complex Hypersurfaces,
{\nineit Annals of Mathematics Studies} {\bf 61}, Princeton University Press,
Princeton, NJ (1968).

\bibitem{NW}
W. Neumann and J. Wahl,
{\nineit  Casson invariant of links of singularities,}
Comment. Math. Helv. {\bf 65} (1990), 58--78.

\bibitem{O}
T. Ohtsuki, {\nineit Ideal points and incompressible surfaces in
two-bridge knot complements,} J. Math. Soc. Japan {\bf 46} (1994),
51--87.

\bibitem{R}
D. Rolfsen,
Knots and Links,
{\nineit Mathematics Lecture Series} {\bf 7},
Publish or Perish, Berkeley, CA (1976).

\bibitem{S}
C. Simpson, {\nineit  Harmonic bundles on noncompact curves,} J.
Amer. Math. Soc. {\bf 3} (1990), 713--770.

\end{thebibliography}
\end{document}